\documentclass[11pt]{article}
\usepackage[centertags]{amsmath}
\usepackage{amsmath}
\usepackage[colorlinks=true,linkcolor=blue,citecolor=blue,urlcolor=blue]{hyperref}
\usepackage{amssymb}
\usepackage{amsthm}
\usepackage{color}
\numberwithin{equation}{section}
\newcommand{\be}{\begin{equation}}
\newcommand{\ee}{\end{equation}}
\newcommand{\dint}{\displaystyle\int}
\newcommand{\1}{1\!\!1}
\newtheorem{definition}{Definition}[section]
\newtheorem{theorem}{Theorem}[section]
\newtheorem{proposition}{Proposition}[section]
\newtheorem{lemma}{Lemma}[section]
\newtheorem{remark}{Remark}[section]

\newcommand{\ind}{\1}
\newcommand{\sgn}{\mbox{sgn}}

\def \R{\mathbb{R}}

\def \N{\mathbb{N}}
\def \E{\mathbb{E}}

\def \L{\mathbb{L}}

\def \bf{\textbf}
\def \it{\textit}
\def \sc{\textsc}
\def \bop {\noindent\textbf{Proof }}

\def \F{\mathcal{F}}

\def \bop {\noindent\textbf{Proof }}

\usepackage{amssymb}
\usepackage{amsthm}
\usepackage{amsmath}
\usepackage{latexsym}
\usepackage[T1]{fontenc}
\usepackage[latin1]{inputenc}

\def \R{\mathbb{R}}

\def \N{\mathbb{N}}
\def \E{\mathbb{E}}

\def \L{\mathbb{L}}

\def \bf{\textbf}
\def \it{\textit}
\def \sc{\textsc}
\def \bop {\noindent\textbf{Proof.}}

\def \F{\mathcal{F}}

\numberwithin{equation}{section}

\def \R{\mathbb{R}}

\def \N{\mathbb{N}}
\def \E{\mathbb{E}}

\def \L{\mathbb{L}}

\def \bf{\textbf}
\def \it{\textit}
\def \sc{\textsc}
\def \bop {\noindent\textbf{Proof }}

\def \F{\mathcal{F}}

\begin{document}

\title{Reflected BSDEs with Logarithmic Growth and Applications\\ in Mixed Stochastic Control Problems \thanks{ This research was supported in part by the ''Scientific Research Projects of Ibn Zohr University''}}
\author{Brahim EL ASRI \thanks{Universit\'e Ibn Zohr, Equipe. Aide \`a la decision,
ENSA, B.P.  1136, Agadir, Maroc. e-mail: b.elasri@uiz.ac.ma }\,\,\,
\, and \, Khalid OUFDIL \thanks{Universit\'e Ibn Zohr, Equipe. Aide \`a la decision,
ENSA, B.P.  1136, Agadir, Maroc. e-mail: khalid.ofdil@gmail.com }\thanks{K. Oufdil is supported in part by the National Center for Scientific and Technical Research (CNRST), Morocco.}}


\date{}
\maketitle \noindent {\bf{Abstract.}} In this article we  study the existence and the uniqueness of a solution for reflected backward stochastic differential
equations in the case when the generator is logarithmic growth in the $z$-variable $(|z|\sqrt{|\ln(|z|)|})$, the  terminal value and obstacle are  an $L^p$-integrable, for a suitable
$p > 2$. To construct the solution we use localization method. We also apply these results to get the existence of an optimal control strategy for the  mixed stochastic control problem in finite horizon. 

\vskip0.2cm
\noindent \text{AMS 2000 Classification subjects:} 60G40, 93E05, 93E20, 60G99.
\vskip0.2cm
\noindent {${Keywords}:$} Reflected backward stochastic differential
equations, Stochastic control, Stochastic differential
games.  

\medskip
\section{Introduction}
\noindent In this paper we study reflected BSDEs with the applications to
stochastic control.\\
 \indent We consider a reflected backward stochastic differential
equation with a generator $\varphi$, a terminal condition
$\xi$ and an obstacle process $(L_t)_{t\leq T}$

\begin{equation}\label{zlogz1}
\left\{\begin{array}{l}
Y_t =\xi + \int_t^T
\varphi(s,Y_s,Z_s) ds +K_T-K_t- \int_t^T Z_s dB_s,\qquad
t\in[0,T]\\
 Y_t \geq L_t,\qquad t\in[0,T]~~~and~~~\int_{0}^{T}\left(Y_s-L_s\right)dK_s=0,
\end{array}\right.
\end{equation}
where $(B_t)_{t\geq 0}$ is a standard
Brownian motion and the process  $K$ is non-decreasing and its role is to push upward $Y$ in order to keep
it above the obstacle $L$. The one barrier reflected BSDEs, have been introduced by
El Karoui et al.  \cite{EKPPQ}.\\

 Reflected BSDEs 
 were studied by several authors (see e.g. \cite{BY2, EKPPQ, HLW, HP}, and the references therein). The motivations are mainly related to applications especially the pricing of American options in markets constrained or not, mixed control, partial
differential variational inequalities, real options, switching optimal (see e.g. \cite{BE,B, BY, BY1, DE, BE1, EH, E, KZ, KZ1}, and the references therein). Once more under square integrability of the data and Lipschitz property of the coefficient $\varphi$, the authors of \cite{EKPPQ} showed the existence and the uniqueness of the
solution.\\

There have been a lot of works which dealt with the issue of existence and uniqueness results under weaker assumptions than the ones of Pardoux and Peng \cite{PP} or El Karoui et al. \cite{EKPPQ}. However, for their own reasons, authors focused only on the weakness of the Lipschitz property of the coefficient and not on the square integrability of the
data $\xi$ and $\varphi$. Actually, there have been relatively few papers which dealt
with the problem of existence and uniqueness of the solution for reflected BSDEs in the case when the
coefficients are not square integrable. Nevertheless, we should point out that Briand et al. (2003, \cite{BDHPS})  have proved existence and uniqueness of a solution for the standard BSDEs  in the case when the data belong only to $L^p$ for some
$p \in ]1, 2[$. Recently,  Hamadene and Popier \cite{HP}  considered reflected BSDEs in the case when the data belong only to $L^p$ for some
$p \in ]1, 2[$. We  proved the existence and the uniqueness of the solution for (\ref{zlogz1}). Recently, Bahlali and El Asri \cite{BE} have considered the BSDE when $\varphi$ is allowed to
have logarithmic growth ($|z|\sqrt{|\ln(|z|)|}$) in the $z$-variable. Moreover, the terminal value is assumed
to be merely $L^p$-integrable, with some $p > 2$.  The main objective of our paper is to study
the existence and the uniqueness of a solution for reflected BSDEs (\ref{zlogz1}) when the generator and the terminal value are as the same as  in  Bahlali and El Asri \cite{BE} and Bahlali et al. \cite{BKK} respectively. These kind of generators are between the linear growth and the quadratic one. In this case, the square integrability of the terminal datum is not sufficient to ensure the existence of solutions while the exponential integrability seems strong enough and is somehow restrictive. In this paper, one should require some $p$-integrability of the terminal datum $\xi$ with $p > 2$. It should be noted that we do not need the comparison theorem in our proofs. The main motivation of this work is its several applications in: finance, control, games, PDEs,....\\

 \bf{ Application in mixed stochastic control problem with finite horizon}\\\\
 Suppose that we have a system, whose evolution is described by the process $X$, which has an
effect on the wealth of a controller. On the other hand the controller has no influence on the
system. The process $X$ may represent, for example, the price of an asset on the market and the
controller a small share holder or a small investor. The controller acts to protect his advantages
by means of $u \in \mathcal{U} $ via the probability $P^u$, here $\mathcal{U}$ is the set of admissible controls. On the other
hand he has also the possibility at any time $\tau \in \mathcal{T}$ to stop controlling. The control is not free.
We define the payoff
  $$ J(u,\tau)=\E^u\left[\int_{0}^{\tau}h(s,X,u_s)ds+g(\tau,X)\ind_{\{\tau<T\}}+g_1(X_{\tau})\ind_{\{\tau=T\}}\right],$$
   $h(.,X,u)$ is the instantaneous reward for the controller,  $g(.,X)$ and $g_1(X)$ are
respectively the rewards if he decides to stop before or until finite time $T$. The problem is to look
for an optimal strategy for the controller i.e. a strategy $(u^*,\tau^*)$ such that
$$J(u,\tau)\leq J(u^*,\tau^*),\qquad \forall(u,\tau)\in \mathcal{U}\times\mathcal{T}.$$
The Hamiltonian associated with this mixed
stochastic control problem is
$$ H(t,X,z,u_t):= z\sigma^{-1}(t,X)f(t,X,u_t)+h(t,X,u_t)\quad \forall(t,X,z,u_t)\in [0,T]\times\Omega\times\R^{d}\times
\mathcal{A}.$$
where $f(t,X,u_t)$ is the drift of dynamics $(X_t)_{t\leq T}$ and $\mathcal{A}$  a compact metric space. The Hamiltonian function attains its
supremum over the set $\mathcal{A}$ at some $u^* \equiv u^*(t,X,z)\in \mathcal{A}$, for
any given $(t,X,z)\in[0,T]\times \Omega \times \R^d$, namely,
$$
\sup\limits_{u\in \mathcal{A}}H(t,X,z,u)=H(t,X,z,u^*(t,X,z)).
$$
The main objective of mixed stochastic control problem is to show the existence of an optimal strategy for
the stochastic control of diffusion. The main idea consists to characterize the value function as
the unique solution of this reflected BSDE
$$
\left\{\begin{array}{l}
Y_t =g_1(X_T) + \int_t^T
H(s,X_s,Z_s,u^*(s,X_s,Z_s)) ds +K_T-K_t- \int_t^T Z_s dB_s,\qquad
t\in[0,T]\\
 Y_t\geq g(t,X_t) ,\qquad \forall t\in[0,T]~~~and~~~\int_{0}^{T}\left(Y_s-g(s,X_s)\right)dK_s=0.
\end{array}\right.
$$
 This problem has been previously studied by Karatzas and Zamfirescu \cite{KZ, KZ1} with $f$, $g_1$, $g$ and $h$ are bounded, who applied the martingale methods. Recently, Bayraktar and Yao \cite{BY1} extended the
work of Karatzas and Zamfirescu \cite{KZ}   but allowing the functions $g_1$, $g$ and $h$ to be bounded below.
Our aim in this paper is to relax the boundedness assumption and replace it by the linear-growth condition.\\

The paper is organized as follows: In section 2, we present the assumptions, we formulate the problem and give some examples. In section 3, we show priori estimates for solutions of reflected BSDEs.
 In section 4, we give estimate between two solutions. In section 5, we give
the main result on the existence and the uniqueness of the solution of reflected BSDEs. Finally, in section 6, we
introduce the mixed stochastic control problem  and we give the
connection between mixed stochastic control problem and reflected BSDEs, we show the value function as a solution of reflected BSDEs. 


\section{Assumptions, setting of the problem and examples}


\subsection{Assumptions}
Let $(\Omega, \mathcal{F}, P)$ be a fixed probability space on which
is defined a standard $d$-dimensional Brownian motion
$B=(B_t)_{0\leq t\leq T}$ whose natural filtration is
$(\mathcal{F}_t^0:=\sigma \{B_s, s\leq t\})_{0\leq t\leq T}$. Let
$\mathbf{F}=(\mathcal{F}_t)_{0\leq t\leq T}$ be
the completed filtration of $(\mathcal{F}_t^0)_{0\leq t\leq T}$ with the $P$%
-null sets of $\mathcal{F}$. We consider the following assumptions,

\begin{enumerate}
\item[(\bf {H.1})] \ \ There exists a positive constant $\lambda$ large enough such that \ \ $\E \left[ |\xi|^{e^{\lambda T}+1}\right]<+\infty.$

\item[(\bf {H.2})]
\begin{itemize}
\item[(i)] The process $(L_t)_{t\leq T}$ is continuous,
\item[(ii)] $\forall p\in]1,2[$, $\E\left[\sup\limits_{0\leq t\leq T}\left((L^+_t)^{e^{\lambda t}}\right)^{\frac{p}{p-1}}\right]<+\infty$.
\end{itemize}

\item[(\bf{H.3})]
\begin{enumerate}\item[(i)]  $\varphi$ is continuous in $(y,z)$ for almost all $(t,\omega)$.
    \item[(ii)] There exists a process $(\eta_t)_{0\leq t\leq T} $ satisfying $$\E\left[\int_{0}^{T} |\eta_s|^{e^{\lambda s}+1}ds\right]<+\infty.$$
    \item[(iii)] There exists a positive constant  $c_0$ such that: for every $t, \omega, y, z$,
    $$\mid \varphi(t,\omega,y,z)\mid \leq |\eta_t|+c_0|z|\sqrt{\vert\ln(|z|)\vert}. $$
    \end{enumerate}
\item[(\bf {H.4})]
There exist $ v\in\L^{q'}(\Omega\times [0, T]; \R_
 +)$ (for some $q'>0$), a real valued sequence $(A_N)_{N>1}$ and constants
 $M_2 \in\R_+$, $r>0$ such that:
\begin{enumerate}
\item[(i)] $\forall N>1$, \quad $1<A_N\leq N^{r}.$
\item[(ii)] $\lim\limits_{N\rightarrow+\infty} A_N =
+\infty .$
\item[(iii)] For every $N\in\N,\; \hbox {and every} \ y,\; y'\; z,\; z'
\;\hbox{such that}\; \mid y\mid,\; \mid y'\mid,\; \mid z\mid, \;\mid
z'\mid\leq N$, we have:
\begin{align*}
\big(y-y^{\prime}\big)
\big(\varphi(t,\omega,y,z)-\varphi(t,\omega,y^{\prime},z')\big)
\1_{\{v_t(\omega)\leq N\}}
 & \leq \
  M_2\mid
y-y^{\prime}\mid^{2}\ln A_{N} \\
& + M_2\mid y-y^{\prime}\mid\mid
z-z^{\prime}\mid\sqrt{\ln A_{N}} \\
& + M_2 \dfrac{\ln A_{N}}{ A_{N}}.
\end{align*}
\end{enumerate}
\end{enumerate}
Next we define the following spaces: for $p>1$,
\begin{itemize}
\item $\mathcal{S}^p$ is the space of $\R $-valued $\mathcal{F}_t$-adapted and continuous processes $\left(Y_t\right)_{t\in[0,T]}$ such that $$||Y||_{\mathcal{S}^p}=\E\left[\sup_{0\leq t\leq T}|Y_t|^p\right]^{\frac{1}{p}}<+\infty.$$
\item $\mathcal{M}^p$ denote the set $\mathcal{P}$-measurable processes $\left(Z_t\right)_{t\in[0,T]}$ with value in $\R^d$  such that:
$$||Z||_{\mathcal{M}^p}=\E\left[\left(\int_{0}^{T}|Z_s|^2ds\right)^{\frac{p}{2}}\right]^{\frac{1}{p}}<+\infty.$$
\item $\mathcal{H}^p$ be the set of adapted continuous non decreasing processes $\left(K_t\right)_{t\in[0,T]}$ such that $K_0=0$ and $\E\left[(K_T)^p\right]<+\infty$.
\end{itemize}
Now it will be convenient to define the notion of solution of the reflected BSDE associated with the triple $(\xi, \varphi, L)$  which we consider throughout this paper.



\begin{definition}\label{def}
We say that $(Y_t,Z_t,K_t)_{\{0\leq t\leq T\}}$ is a solution of the reflected BSDE associated with the terminal condition
$\xi$, the generator $\varphi$ and the obstacle $L$ if the followings hold:
\begin{enumerate}
\item[1-]$(Y,Z,K)\in\mathcal{S}^{e^{\lambda T}+1}\times\mathcal{M}^2\times \mathcal{H}^p$, $\forall p\in]1,2[$;
\item[2-]$Y_t = \xi + \int_t^T
\varphi(s,Y_s,Z_s) ds +K_T-K_t- \int_t^T Z_s dB_s,\qquad t\in[0,T];$
\item[3-]$L_t\leq Y_t, \qquad t\in[0,T]$;
\item[4-]$\int_0^T(Y_s-L_s)dK_s=0 $.
\end{enumerate}
\end{definition}

Now we shall present some examples mentioned in \cite{BEHE} of functions $\varphi$ that satisfies our assumptions.
\subsection{Examples}
\textbf{Example 1.} Let $\varphi\in{\cal C}(\R^{d};\R_+)\cap{\cal C}^1(\R^{d}-\{0\};\R_+)$ be such that: for $\varepsilon\in]0,1[,$

\begin{equation}
\varphi(z)=\left\{\begin{array}{l}
|z|\sqrt{-\ln(|z|)}\qquad\qquad \text{ if } |z|<1-\varepsilon
\\
|z|\sqrt{\ln(|z|)}\qquad\qquad\quad \text{ if } |z|>1+\varepsilon.
\end{array}\right.
\end{equation}
It is not difficult to see that $\varphi$ satisfies (\textbf{H.3}).\\

To verify that $\varphi$ satisfies also (\textbf{H.4}), it is enough to show that for every $z$ and $z'$ such that $|z|,|z'|\leq N$
\begin{equation}
|\varphi(z)-\varphi(z')|\leq c\left(\sqrt{\ln (N)}|z-z'|+\frac{\ln (N)}{N}\right)
\end{equation}
for $N$ large enough and some positive constant $c$. This can be proved by considering separately the following four cases. $0\leq |z|,|z'|\leq \frac{1}{N}$, $\frac{1}{N}\leq |z|,|z'|\leq 1-\varepsilon$, $1-\varepsilon\leq|z|,|z'|\leq 1+\varepsilon$ and $1+\varepsilon\leq |z|,|z'|\leq N$.\\
In the first case for $N>\sqrt{e}$, the map $x\rightarrow x\sqrt{-\ln x}$ increases for $x\in]0,\frac{1}{N}]$, it follows that
$$|\varphi(z)-\varphi(z')|\leq |\varphi(z)|+|\varphi(z')|\leq 2\frac{1}{N}\ln (N).$$
To prove the result for the rest of the cases, we apply the finite increment theorem on $\varphi$.\\

The following example shows that our assumptions enable to treat some reflected BSDEs with stochastic monotone coefficient.\\

\noindent
\textbf{Example 2.} Let $\varphi$ satisfying (\textbf{H.3}) and
\begin{equation*}
(\textbf{H'.4})\left\{\begin{array}{l}
\text{There exist a positive process } (C_t)_{0\leq t\leq T} \text{ satisfying }\E\left[\int_{0}^{T}e^{q'C_s}ds\right]<+\infty  \\ (\text{for some }q'>0) \text{ and } M\in\R_+ \text{ such that:}
\\  (y-y')(\varphi(t,\omega,y,z)-\varphi(t,\omega,y',z')) \leq M|y-y'|^2(C_t(\omega)+|\ln |y-y'||)\\ \qquad\qquad\qquad\qquad\qquad\qquad\qquad\quad+M|y-y'||z-z'|\sqrt{C_t(\omega)+|\ln|z-z'||}.
\end{array}\right.
\end{equation*}
In particular we have for all $z,z'$
\begin{equation*}
|\varphi(t,\omega,y,z)-\varphi(t,\omega,y,z')|\leq M|z-z'|\sqrt{C_t(\omega)+|\ln|z-z'||}.
\end{equation*}
To see if (\textbf{H.4}) is verified, we consider the following cases
$$|y-y'|\leq \frac{1}{2N},\qquad \frac{1}{2N}\leq |y-y'|\leq 2N,$$
and
$$|z-z'|\leq \frac{1}{2N},\qquad \frac{1}{2N}\leq |z-z'|\leq 2N.$$
It is not difficult to prove that for some constant $c$ we have:
\begin{eqnarray*}
&&  (y-y')(\varphi(t,y,z)-\varphi(t,y',z)) \leq c\ln (N)\left(|y-y'|^2+\frac{1}{N}\right)\\
&& |\varphi(t,y,z)-\varphi(t,y,z')|\leq c\sqrt{\ln (N)}\left(|z-z'|+\frac{1}{N}\right),
\end{eqnarray*}
whenever $v_s:=e^{C_s}\leq N$ and $|y|,~|y'|,~|z|,~|z'|\leq N$.\\

\section{Apriori estimates}

Here, we want to obtain estimates for solutions to  reflected BSDEs  in the spirit of the work
\cite{BKK} which shows that these estimates are very useful in the study of existence and uniqueness of
solutions. 
 Before we do so, we give the following lemma that has already been mentioned and proved in \cite{BKK}, but for the sake of the reader we give it again.
\begin{lemma}\label{important}
Let $(y,z)\in\R\times \R^d$ be such that $y$ is large enough. Then, for every $C_1>0$ there exists $C_2>0$ such that:
\begin{equation}\label{yleqz}
C_1\mid y\mid|z|\sqrt{\vert\ln(|z|)\vert}\leq \frac{|z|^2}{2}+ C_2\ln(\mid
y\mid)\mid y\mid^{2}.
\end{equation}
\end{lemma}

Now we start by showing how  to control the process $Y$ in terms of the data and $K_T$.
\begin{lemma} \label{estimateY}
Let $(Y,Z,K)$ be a solution of the reflected BSDE \eqref{zlogz1}, where $(\xi,\varphi,L)$
satisfies the assumptions $(\bf {H.1})$, $(\bf {H.2})$ and $(\bf {H.3})$. Then, for $\lambda>1$, there
exists a constant $C(\lambda,T)$ such that:
\begin{equation}\label{estY}
\E \left[\sup_{t \in [0,T]} |Y_t|^{e^{\lambda t}+1}\right] \leq C(\lambda,T)\E \left[ |\xi|^{e^{\lambda T}+1} + \int_0^{T}|\eta_s|^{e^{\lambda s}+1}ds+
\sup_{t \in [0,T]}\left(
(L_t^+)^{e^{\lambda t}}\right)K_T \right ].
\end{equation}
\end{lemma}

\noindent
\textbf{Proof.} Without loss of generality we assume that the $y$-variable is sufficiently large. For some constant $\lambda$ large enough, let us consider the function from
$[0,T]\times\R$ into $\R^+$ defined by
$$u(t,x)=\mid x\mid^{e^{\lambda t}+1}.$$
Then,
\begin{itemize}
\item $u_t=\lambda e^{\lambda t}\ln(\mid x\mid)\mid x\mid^{e^{\lambda t}+1}$;
\item $ u_x=(e^{\lambda t}+1)\mid x\mid^{e^{\lambda t}}\sgn(x)$;
\item $u_{xx}=(e^{\lambda t}+1)e^{\lambda t}\mid x\mid^{e^{\lambda t}-1}$;
\end{itemize}
where $\sgn(x)=-\ind_{x\leq0}+\ind_{x>0}$.\\

\noindent
Let $(\tau_{k})_{k\geq 0}$ be the sequence of stopping times defined as follows:
$$\tau_k:=\inf \left\{t\geq 0,\,\left[ \int_{0}^{t}(e^{\lambda s}+1)^2\mid
Y_s\mid^{2e^{\lambda s}}\mid Z_s\mid^2 ds\right]\vee|Y_t|\geq k\right\}\wedge T.$$
Now we apply It\^o's formula on the process $Y$ and the function $y\mapsto \mid y\mid^{e^{\lambda t}+1}$ to obtain :
\begin{eqnarray*}
&& \mid Y_{t\wedge \tau_k}\mid^{e^{\lambda (t\wedge\tau_k)}+1} =\mid
Y_{\tau_k}\mid^{1+e^{\lambda \tau_k}}- \lambda\int_{t\wedge \tau_k}^{\tau_k}e^{\lambda s}\ln(\mid Y_s \mid)\mid  Y_s\mid^{e^{\lambda s}+1}ds \\
\nonumber
&& \qquad\qquad \qquad~~~ -\frac{1}{2}\int_{t\wedge
\tau_k}^{\tau_k} |Z_s|^2(e^{\lambda s}+1)e^{\lambda s}\mid
Y_s\mid^{e^{\lambda s}-1}ds
\\\nonumber
&& \qquad \qquad\qquad~~~ +\int_{t\wedge \tau_k}^{\tau_k}
(e^{\lambda s}+1)\mid Y_s\mid^{e^{\lambda s}}\sgn(Y_s)\varphi(s,Y_s,Z_s)ds  \\
\\\nonumber
&& \qquad \qquad\qquad~~~ +\int_{t\wedge \tau_k}^{\tau_k}
(e^{\lambda s}+1)\mid Y_s\mid^{e^{\lambda s}}\sgn(Y_s)dK_s  \\
\nonumber
&& \qquad \qquad\qquad~~~ -  \int_{t\wedge
\tau_k}^{\tau_k}(e^{\lambda s}+1)\mid
Y_s\mid^{e^{\lambda s}}\sgn(Y_s) Z_s dB_s.
\end{eqnarray*}
Using (iii) in assumption (\bf {H.3}) we have:
\begin{eqnarray*}
&& \mid Y_{t\wedge \tau_k}\mid^{e^{\lambda (t\wedge\tau_k)}+1} \leq \mid
Y_{\tau_k}\mid^{1+e^{\lambda \tau_k}}- \lambda \int_{t\wedge \tau_k}^{\tau_k}e^{\lambda s}\ln(\mid Y_s \mid)\mid  Y_s\mid^{e^{\lambda s}+1}ds \\
\nonumber
&& \qquad \qquad\qquad~~~ -\frac{1}{2}\int_{t\wedge
\tau_k}^{\tau_k} |Z_s|^2(e^{\lambda s}+1)e^{\lambda s}\mid
Y_s\mid^{e^{\lambda s}-1}ds
\\\nonumber
&& \qquad \qquad\qquad~~~ +\int_{t\wedge \tau_k}^{\tau_k}
(e^{\lambda s}+1)\mid Y_s\mid^{e^{\lambda s}}(|\eta_s|+c_0|Z_s|\sqrt{\vert\ln(|Z_s|)\vert})ds  \\
\\\nonumber
&& \qquad \qquad\qquad~~~ +\int_{t\wedge \tau_k}^{\tau_k}
(e^{\lambda s}+1)\mid Y_s\mid^{e^{\lambda s}}\sgn(Y_s)dK_s  \\
\nonumber
&& \qquad \qquad\qquad~~~ -  \int_{t\wedge
\tau_k}^{\tau_k}(e^{\lambda s}+1)\mid Y_s\mid^{e^{\lambda s}}\sgn(Y_s) Z_s
dB_s.
\end{eqnarray*}
Next by Young's inequality we have:
$$
(e^{\lambda s}+1)\mid Y_s\mid^{e^{\lambda s}}|\eta_s| \leq \mid
Y_s\mid^{e^{\lambda s}+1}+ (e^{\lambda s}+1)^{e^{\lambda s}+1}|\eta_s|^{e^{\lambda s}+1}.
$$
For $\mid Y_s \mid$ large enough and thanks to the last inequality, it follows that:
\begin{eqnarray*}
&& \mid Y_{t\wedge \tau_k}\mid^{e^{\lambda (t\wedge \tau_k)}+1} \leq \mid
Y_{\tau_k}\mid^{1+e^{\lambda \tau_k}}-\lambda \int_{t\wedge \tau_k}^{\tau_k}e^{\lambda s}\ln(\mid Y_s \mid)\mid  Y_s\mid^{e^{\lambda s}+1}ds \\
\nonumber && \qquad\qquad\qquad\quad  -\frac{1}{2}\int_{t\wedge
\tau_k}^{\tau_k} |Z_s|^2(e^{\lambda s}+1)e^{\lambda s}\mid
Y_s\mid^{e^{\lambda s}-1}ds
\\\nonumber && \qquad\qquad\qquad\quad   +\int_{t\wedge \tau_k}^{\tau_k}
c_0(e^{\lambda s}+1)\mid Y_s\mid^{e^{\lambda s}}|Z_s|\sqrt{\vert\ln(|Z_s|)\vert}ds
\\\nonumber && \qquad\qquad\qquad\quad
+\int_{t\wedge\tau_k}^{\tau_k}|Y_s|^{e^{\lambda s}+1}\ln(|Y_s|)ds +\int_{t\wedge \tau_k}^{\tau_k}(e^{\lambda s}+1)^{e^{\lambda s}+1}|\eta_s|^{ e^{\lambda s}+1}ds \\
 \\\nonumber && \qquad\qquad\qquad\quad   +\int_{t\wedge \tau_k}^{\tau_k}
(e^{\lambda s}+1)\mid Y_s\mid^{e^{\lambda s}}\sgn(Y_s)dK_s  \\
\nonumber && \qquad\qquad\qquad\quad   -  \int_{t\wedge
\tau_k}^{\tau_k}(e^{\lambda s}+1)\mid
Y_s\mid^{e^{\lambda s}}\sgn(Y_s) Z_s dB_s.
\end{eqnarray*}
Then we have:
\begin{eqnarray*}
&& \mid Y_{t\wedge \tau_k}\mid^{e^{\lambda (t\wedge \tau_k)}+1} \leq \mid
Y_{\tau_k}\mid^{1+e^{\lambda \tau_k}}+ \int_{t\wedge \tau_k}^{\tau_k}|Y_s|^{e^{\lambda s}-1}\left[c_0(e^{\lambda s}+1)\mid Y_s\mid|Z_s|\sqrt{\vert\ln(|Z_s|)\vert}\right.\\
&& \qquad\qquad\qquad\quad \left.-(e^{\lambda s}+1)e^{\lambda s}\frac{|Z_s|^2}{2}-(\lambda e^{\lambda s}-1)\ln(\mid Y_s \mid)\mid  Y_s\mid^2\right]ds \\
\\\nonumber && \qquad\qquad\qquad\quad
 +\int_{t\wedge \tau_k}^{\tau_k}(e^{\lambda s}+1)^{e^{\lambda s}+1}|\eta_s|^{ e^{\lambda s}+1}ds \\
 \\\nonumber && \qquad\qquad\qquad\quad   +\int_{t\wedge \tau_k}^{\tau_k}
(e^{\lambda s}+1)\mid Y_s\mid^{e^{\lambda s}}\sgn(Y_s)dK_s  \\
\nonumber && \qquad\qquad\qquad\quad   -  \int_{t\wedge
\tau_k}^{\tau_k}(e^{\lambda s}+1)\mid
Y_s\mid^{e^{\lambda s}}\sgn(Y_s) Z_s dB_s.
\end{eqnarray*}
Since  $\lambda e^{\lambda s}-1>0$ for  $\lambda>1$, then using Lemma \ref{important}, we have for $\lambda$ large enough:
\begin{equation}
c_0(e^{\lambda s}+1)\mid Y_s\mid|Z_s|\sqrt{\vert\ln(|Z_s|)\vert}\leq(e^{\lambda s}+1)e^{\lambda s}\frac{|Z_s|^2}{2}+(\lambda e^{\lambda s}-1)\ln(\mid Y_s \mid)\mid  Y_s\mid^2.
\end{equation}
Therefore, it follows that:
\begin{eqnarray}\label{key}
&& \mid Y_{t\wedge \tau_k}\mid^{e^{\lambda (t\wedge \tau_k)}+1} \leq \mid
Y_{\tau_k}\mid^{1+e^{\lambda \tau_k}}+\int_{t\wedge \tau_k}^{\tau_k}(e^{\lambda s}+1)^{e^{\lambda s}+1}|\eta_s|^{ e^{\lambda s}+1}ds
 \\\nonumber && \qquad\qquad\qquad\quad   +\int_{t\wedge \tau_k}^{\tau_k}
(e^{\lambda s}+1)\mid Y_s\mid^{e^{\lambda s}}\sgn(Y_s)dK_s  \\
\nonumber && \qquad\qquad\qquad\quad   -  \int_{t\wedge
\tau_k}^{\tau_k}(e^{\lambda s}+1)\mid
Y_s\mid^{e^{\lambda s}}\sgn(Y_s) Z_s dB_s.
\end{eqnarray}
Next let us deal with the term $\int_{t\wedge \tau_k}^{\tau_k}
(e^{\lambda s}+1)\mid Y_s\mid^{e^{\lambda s}}\sgn(Y_s)dK_s $.
Indeed, the hypothesis related to increments of $K$ and $Y-L$ implies that:
$dK_s=\ind_{\{Y_s = L_s\}}dK_s$, for any $s\leq T$. Therefore we have:
\begin{eqnarray*}
&& \int_{t\wedge \tau_k}^{\tau_k} (e^{\lambda s}+1)\mid
Y_s\mid^{e^{\lambda s}}\sgn(Y_s)dK_s
 =\int_{t\wedge \tau_k}^{\tau_k}
(e^{\lambda s}+1)\mid Y_s\mid^{e^{\lambda s}-1}Y_s \ind_{\{Y_s =
L_s\}}dK_s\\\nonumber
&& \qquad\qquad\qquad\qquad\quad \qquad  \qquad ~~~~~~ =\int_{t\wedge \tau_k}^{\tau_k}
(e^{\lambda s}+1)\mid L_s\mid^{e^{\lambda s}-1}L_s \ind_{\{Y_s =
	L_s\}}dK_s.
\end{eqnarray*}
It follows that:
\begin{eqnarray*}
&& \int_{t\wedge \tau_k}^{\tau_k} (e^{\lambda s}+1)\mid
Y_s\mid^{e^{\lambda s}}\sgn(Y_s)dK_s 	
\leq \int_{t\wedge \tau_k}^{\tau_k}
(e^{\lambda s}+1)(L_s^+)^{e^{\lambda s}} dK_s
\\\nonumber && \qquad \qquad  \qquad \qquad\qquad\qquad\qquad~ \leq 2(e^{\lambda T}+1)
\sup_{t \in [0,T]}\left(
(L_t^+)^{e^{\lambda t}}\right)K_{\tau_k}\\\nonumber &&  \qquad \qquad  \qquad \qquad\qquad\qquad\qquad~ \leq 2(e^{\lambda T}+1)
\sup_{t \in [0,T]}\left(
(L_t^+)^{e^{\lambda t}}\right)K_{T},
\end{eqnarray*}
which means that:
\begin{equation}\label{Kt}
\int_{t\wedge \tau_k}^{\tau_k} (e^{\lambda s}+1)\mid
Y_s\mid^{e^{\lambda s}}\sgn(Y_s)dK_s\leq 2(e^{\lambda T}+1)
\sup_{t \in [0,T]}\left(
(L_t^+)^{e^{\lambda t}}\right)K_{T}.
\end{equation}
We combine (\ref{key}) and (\ref{Kt}) and we take expectation to get:
\begin{eqnarray}
 &&\E\left[\mid Y_{t\wedge \tau_k}\mid^{e^{\lambda (t\wedge \tau_k)}+1}\right] \leq \E \left[\mid
Y_{\tau_k}\mid^{1+e^{\lambda \tau_k}}\right]+(e^{\lambda T}+1)^{e^{\lambda T}+1}\E \left[\int_{t\wedge \tau_k}^{\tau_k}|\eta_s|^{ e^{\lambda s}+1}ds\right]\\\nonumber
&& \qquad\qquad\qquad\qquad\quad+2(e^{\lambda T}+1)\E \left[
\sup_{t \in [0,T]}\left(
(L_t^+)^{e^{\lambda t}}\right)K_{T}\right].
\end{eqnarray}
Since the sequence of stopping times $(\tau_k)_{k\geq 0}$ is increasing of stationary type, we pass to the limit when $k\rightarrow+\infty$ then we use Fatou's Lemma to obtain:
\begin{eqnarray}
 &&\E\left[\mid Y_{t}\mid^{e^{\lambda t}+1}\right] \leq \E \left[\mid
\xi\mid^{e^{\lambda T}+1}\right]+(e^{\lambda T}+1)^{e^{\lambda T}+1}\E \left[\int_{0}^{T}|\eta_s|^{ e^{\lambda s}+1}ds\right]\\\nonumber
&& \qquad\qquad\qquad~+2(e^{\lambda T}+1)\E \left[
\sup_{t \in [0,T]}\left(
(L_t^+)^{e^{\lambda t}}\right)K_{T}\right].
\end{eqnarray}
Finally we use Burkholder-Davis-Gundy's inequality to complete the proof.\qed\\

We will now establish an estimate for the process $Z$. Actually we have:
\begin{lemma} \label{estimateZ}
Let $(Y,Z,K)$ be a solution of the reflected BSDE \eqref{zlogz1}, where $(\xi,\varphi,L)$
satisfies the assumptions (\bf {H.1}), (\bf {H.2}) and (\bf {H.3}). Then, there exists a positive constant $C(\lambda,c_0,T)$ such that:
\begin{equation}\label{estZ}
\E \left[ \int_0^{T} |Z_s|^2 ds \right]
\leq C(\lambda,c_0,T) \E \left[1+ |\xi|^{e^{\lambda T}+1} + \int_0^T\mid
\eta_s \mid^{e^{\lambda s}+1}ds +\sup_{t \in [0,T]} \left((L^+_t)^{e^{\lambda t}}\right) K_T\right].
\end{equation}
\end{lemma}


\noindent
\textbf{Proof.} Applying It\^o's formula to the process $Y$ and the
function $y \mapsto y^2$ yields:
\begin{eqnarray}\label{estimz}
\nonumber
&& |Y_0|^2 + \int_0^{T} |Z_s|^2 ds = \xi^2 + 2\int_0^{T} Y_s
\varphi(s,Y_s,Z_s) ds + 2\int_0^{T} Y_s dK_s- 2 \int_0^{T}  Y_s Z_s dB_s \\ \nonumber
&& \qquad\qquad\qquad\qquad~\leq |\xi|^2 + 2\int_0^{T}\mid Y_s\mid(\mid
\eta_s\mid+c_0|Z_s|\sqrt{\vert\ln(|Z_s|)\vert}) ds+ 2\int_0^{T}
Y_s  dK_s \\
&&\qquad\qquad\qquad\qquad~ -2 \int_0^{T} Y_s Z_s dB_s .
\end{eqnarray}
Using Lemma \ref{important} we can show that there exists a constant $\tilde{C}$ that depends on  $c_0$ such that
\begin{equation}\label{eps}
2c_0|Y_s||Z_s|\sqrt{\vert\ln(|Z_s|)\vert}\leq \frac{|Z_s|^2}{2}+\tilde{C}\ln(|Y_s|)|Y_s|^2.
\end{equation}
From Young's inequality we have:
\begin{equation}\label{you}
2\mid Y_s\mid \mid \eta_s\mid\leq \mid Y_s\mid^2+\mid
\eta_s\mid^2.
\end{equation}
Then by combining \eqref{eps} and \eqref{you} with \eqref{estimz} and by using the fact that
$\int_0^T(Y_s-L_s)dK_s =0$ we obtain:
\begin{eqnarray*}
&& |Y_0|^2 +\frac{1}{2} \int_0^{T} |Z_s|^2 ds \leq |\xi|^2 + \sup_{0\leq s\leq
T}|Y_s|^2+\int_0^T\mid \eta_s
\mid^2ds+\tilde{C}\int_0^T \ln(|Y_s|)|Y_s|^2 ds\\
&&\qquad\qquad\qquad\qquad~~- 2 \int_0^{T} Y_s Z_s dB_s+2 \int_{0}^{T}L_s dK_s.
\end{eqnarray*}
Since we assume that $|Y_s|$ is large enough, we have for any $\varepsilon>0$:
$$\mid Y_s\mid^2|\ln(|Y_s|)|\leq \mid Y_s\mid^{2+\varepsilon} \mbox{ and } \mid Y_s\mid^{2}\leq \mid Y_s\mid^{2+\varepsilon}.$$
Therefore, there exists a positive constant $\tilde{C}_1$ that depends on $T$ and $c_0$ such that:
\begin{eqnarray}\label{once}
&&\frac{1}{2}\int_0^{T} |Z_s|^2 ds \leq
|\xi|^2 +\tilde{C}_1\sup\limits_{0\leq s\leq T}\mid Y_s\mid^{2+\varepsilon}+\int_0^T\mid \eta_s
\mid^2ds + 2 \int_{0}^{T}L_sdK_s\\ \nonumber
&&\qquad\qquad\qquad- 2\int_0^{T} Y_s Z_s dB_s .
\end{eqnarray}
We choose $\varepsilon=e^{\lambda T}-1$ and we take expectation in both sides of \eqref{once}. Then, there exists a constant $\tilde{C}_2>0$ that still depends on $c_0$ and $T$ such that:
\begin{eqnarray}\label{once1}
\nonumber
&& \E\left[ \int_0^{T} |Z_s|^2 ds\right] \leq \tilde{C}_2\E\left[|\xi|^2 +
\sup\limits_{0\leq s\leq T}\mid Y_s\mid ^{e^{\lambda T}+1} +
\int_0^T\mid \eta_s \mid^2ds+  \int_{0}^{T}L_sdK_s\right]\\
&&\quad\qquad\qquad\qquad+ 2\E\left[\left|\int_0^{T} Y_s Z_s dB_s\right|\right].
\end{eqnarray}
Now, thanks to Burkholder-Davis-Gundy's inequality we have for any $\beta > 0$
\begin{eqnarray*}
&&\E\left[\sup\limits_{t\in[0,T]}\left|\int_t^{T} Y_s Z_s dB_s\right|\right] \leq  \tilde{K} \E\left[ \left(\int_0^{T} |Y_s|^2 |Z_s|^2 ds\right)^{\frac{1}{2}}\right]\\
&&\qquad\qquad\qquad\qquad\qquad~~~\leq \tilde{K} \E\left[\sup\limits_{0\leq s\leq T}\mid Y_s\mid \left(\int_0^{T} |Z_s|^2 ds\right)^{\frac{1}{2}}\right]\\
&&\qquad\qquad\qquad\qquad\qquad~~~\leq \frac{\tilde{K}}{2\beta} \E\left[\sup\limits_{0\leq s\leq T}\mid Y_s\mid^2\right]+ \frac{\beta \tilde{K}}{2} \E\left[\int_0^{T} |Z_s|^2 ds\right].
\end{eqnarray*}
Therefore, choosing $\beta$ small enough, yields to the existence of a positive constant $\tilde{C}_3$ that depends on $c_0$ and $T$ such that:
\begin{eqnarray*}
&&\E \left[\int_0^{T} |Z_s|^2 ds \right]
\leq \tilde{C}_3 \E \left[ |\xi|^2+\sup_{t \in [0,T]} |Y_t|^{e^{\lambda T}+1} + \int_0^T\mid
\eta_s \mid^2ds +\sup_{t \in [0,T]} (L^+_t) K_T\right].
\end{eqnarray*}
Now by using Inequality (\ref{estY}), Young's inequality and H\"{o}lder's inequality, we conclude that there exists a positive constant $C(\lambda,c_0,T)$ such that:
$$\E \left[ \int_0^{T} |Z_s|^2 ds \right]
\leq C(\lambda,c_0,T) \E \left[ 1+|\xi|^{e^{\lambda T}+1} + \int_0^T\mid
\eta_s \mid^{e^{\lambda s}+1}ds +\sup_{t \in [0,T]} \left((L^+_t)^{e^{\lambda t}}\right) K_T\right].$$
The proof is now complete. \qed\\

Now we focus on the control of the process $K.$
\begin{lemma}\label{estimateKK}
Let $(Y,Z,K)$ be a solution of the reflected BSDE \eqref{zlogz1}, where $(\xi,\varphi,L)$
satisfies the assumptions (\bf {H.1}), (\bf {H.2}) and (\bf {H.3}). Then, there
exists a constant $C(\lambda,p,c_0,T)>0$ such that: $\forall p\in]1,2[$,
\begin{equation}\label{estK}
\E \left[ (K_T)^{p}\right] \leq C(\lambda,p,c_0,T)\E \left[1+|\xi|^{e^{\lambda T}+1} +\int_0^T\mid \eta_s \mid^{e^{\lambda s}+1}ds+ \sup_{0\leq t\leq T}\left((L^+_t)^{e^{\lambda t}}\right)^{\frac{p}{p-1}}\right].
\end{equation}
\end{lemma}

\noindent
\textbf{Proof.}
First we recall that:
\begin{equation}\label{ham}
K_T-K_t= Y_t-\xi-\int_{t}^{T} \varphi(s,Y_s,Z_s) ds+\int_{t}^{T} Z_s dB_s.
\end{equation}
Next using  the predictable dual projection property (see e.g. \cite{DM}) we have: for any $t\leq T$,
\begin{eqnarray*}
&&\E \left[\left(K_T-K_t\right)^p\right]=\E\left[\int_{t}^{T}p\left(K_T-K_s\right)^{p-1}dK_s\right]\\\nonumber
&&\qquad\qquad\qquad\quad=p\E \left[\int_{t}^{T}\E\left[\left(K_T-K_s\right)^{p-1}|\mathcal{F}_s\right]dK_s\right].
\end{eqnarray*}
Since $p\in]1,2[$, then thanks to Jensen's conditional inequality we have:
$$\E \left[\left(K_T-K_t\right)^p\right]\leq p\E \left[\int_{t}^{T}\left[\E\left(K_T-K_s\right)|\mathcal{F}_s\right]^{p-1}dK_s\right]. $$
From \eqref{ham} we obtain:
\begin{eqnarray*}
&&\E \left[\left(K_T-K_t\right)^p\right]\leq p\E \left[\int_{t}^{T}\left[\E\left(Y_s-\xi-\int_{s}^{T}\varphi(u,Y_u,Z_u)du\bigg|\mathcal{F}_s\right)\right]^{p-1}dK_s\right]\\
&&\qquad\qquad\qquad\quad  \leq p\E \left[\int_{t}^{T}\left[\E\left(2\sup\limits_{u\in[s,T]}|Y_u|+\int_{s}^{T}|\varphi(u,Y_u,Z_u)|du\bigg|\mathcal{F}_s\right)\right]^{p-1}dK_s\right]\\
&&\qquad\qquad\qquad\quad  =p\E \left[\int_{t}^{T}\Gamma_s^{p-1}dK_s\right]\\
&&\qquad\qquad\qquad\quad  \leq p\E \left[\left(\sup\limits_{s\in[t,T]}|\Gamma_s|\right)^{p-1}\left(K_T-K_t\right)\right];
\end{eqnarray*}
with
$$\Gamma_s=\E\left(2\sup\limits_{u\in[s,T]}|Y_u|+\int_{s}^{T}|\varphi(u,Y_u,Z_u)|du\bigg|\mathcal{F}_s\right).$$
Once more by using Young's inequality  we have:
\begin{eqnarray*}
&& p\E \left[\left(\sup\limits_{s\in[t,T]}|\Gamma_s|\right)^{p-1}\left(K_T-K_t\right)\right]\\
&&\qquad\qquad\qquad\qquad\leq\frac{1}{2}\E \left[\left(K_T-K_t\right)^p\right]+(p-1)2^{\frac{1}{p-1}}\E\left[\left(\sup\limits_{s\in[t,T]}|\Gamma_s|\right)^p\right].
\end{eqnarray*}
Hence
\begin{eqnarray*}
&& \E \left[\left(K_T-K_t\right)^p\right]\leq \frac{1}{2}\E \left[\left(K_T-K_t\right)^p\right]\\
&& \qquad\qquad\qquad\quad+C_p\E\sup\limits_{s\in[t,T]}\left[\E\left(2\sup\limits_{u\in[s,T]}|Y_u|+\int_{s}^{T}|\varphi(u,Y_u,Z_u)|du\bigg|\mathcal{F}_s\right)\right]^p.
\end{eqnarray*}
Thus, using Doob's maximal inequality we obtain:
\begin{eqnarray*}
&& \frac{1}{2}\E \left[\left(K_T-K_t\right)^p\right]\leq C_p \sup\limits_{s\in[t,T]}\E\left[\E\left(2\sup\limits_{u\in[t,T]}|Y_u|+\int_{t}^{T}|\varphi(u,Y_u,Z_u)|du\bigg|\mathcal{F}_s\right) \right]^p\\
&& \qquad\qquad\qquad\quad~\leq C'_p\E \left[\sup\limits_{u\in[t,T]}|Y_u|^p+\left(\int_{t}^{T}|\varphi(u,Y_u,Z_u)|du\right)^p\right].
\end{eqnarray*}
Then by taking $t=0$, there exists a constant $C''_p>0$ such that:
\begin{equation}
\E \left[(K_T)^p\right]\leq  C''_p\E \left[\sup\limits_{s\in[0,T]}|Y_s|^p+\left(\int_{0}^{T}|\varphi(s,Y_s,Z_s)|ds\right)^p\right].
\end{equation}

Using H\"{o}lder's inequality and then assumption \bf{(H.3)}, we obtain that there exists a positive constant $C(p,T)$ (changes from line to line) such that:
\begin{eqnarray*}
&&\left( \int_{0}^{T} \mid \varphi(s,Y_s,Z_s)\mid ds \right)^{p}\leq C(p,T) \int_{0}^{T} \mid \varphi(s,Y_s,Z_s)\mid  ^{p}ds\\
&&\qquad\qquad\qquad\qquad\qquad~~~~\leq C(p,T)\left(
\int_{0}^{T}|\eta_s|^p ds+\int_{0}^{T}\left( c_0|Z_s|\sqrt{\vert\ln(|Z_s|)\vert} \right)^p ds\right),
\end{eqnarray*}
and for any $\varepsilon>0$ we have:
$$\sqrt{2\varepsilon \vert\ln(|Z_s|)\vert}=\sqrt{\vert\ln(|Z_s|^{2\varepsilon }\vert)}\leq|Z_s|^{\varepsilon }.$$
Therefore,
\begin{equation}\label{varphi}
\left( \int_{0}^{T} \mid \varphi(s,Y_s,Z_s)\mid ds \right)^{p}\leq  C(p,T)\left(\int_{0}^{T}|\eta_s|^p ds+\left(\frac{c_0}{\sqrt{2\varepsilon}}\right)^p\int_{0}^{T}|Z_s|^{p(1+\varepsilon)} ds\right).
\end{equation}
We now put $\varepsilon=\frac{2}{p}-1$, then, there exists  a constant  $C(p,c_0,T)>0$ such that:
 $$ \left( \int_{0}^{T} \mid \varphi(s,Y_s,Z_s)\mid ds \right)^{p}\leq  C(p,c_0,T)\left(
\int_{0}^{T}|\eta_s|^p ds+\int_{0}^{T}|Z_s|^2 ds\right).$$
Hence, there exists a positive constant which we still denote $C(p,c_0,T)$ such that:
\begin{equation}\label{maj_K} \E \left[(K_T)^{p}\right] \leq C(p,c_0,T)\E \left[\sup\limits_{s\in[0,T]}|Y_s|^p+ \int_{0}^{T}|\eta_s|^p ds+\int_{0}^{T}|Z_s|^2ds\right].\end{equation}
Using Inequalities \eqref{estY} and (\ref{estZ}), Young's inequality and H\"{o}lder's inequality, there exists a positive constant $C(\lambda,p,c_0,T)$ that changes from line to line such that:
\begin{eqnarray*}
&& \E \left[(K_T)^{p}\right]\leq C(\lambda,p,c_0,T)\E \left[1+|\xi|^{e^{\lambda T}+1} +\int_0^T\mid \eta_s \mid^{e^{\lambda s}+1}ds+\sup_{t \in [0,T]} \left((L^+_t)^{e^{\lambda t}}\right) K_T\right]\\
&&\qquad\qquad\leq C(\lambda,p,c_0,T)\E \left[1+|\xi|^{e^{\lambda T}+1} +\int_0^T\mid \eta_s \mid^{e^{\lambda s}+1}ds+ \sup_{0\leq t\leq T}\left((L^+_t)^{e^{\lambda t}}\right)^{\frac{p}{p-1}}\right]\\
&&\qquad\qquad+\frac{1}{p}\E\left[(K_T)^p\right].
\end{eqnarray*}
\noindent
Finally, there exists a positive constant which we still denote $C(\lambda,p,c_0,T)$ such that:
$$ \E \left[ (K_T)^{p}\right] \leq C(\lambda,p,c_0,T)\E \left[1+|\xi|^{e^{\lambda T}+1} +\int_0^T\mid \eta_s \mid^{e^{\lambda s}+1}ds+ \sup_{0\leq t\leq T}\left((L^+_t)^{e^{\lambda t}}\right)^{\frac{p}{p-1}}\right].$$
The proof is now complete. \qed
\begin{proposition}\label{estimateK}
Let $(Y,Z,K)$ be a solution of the reflected BSDE \eqref{zlogz1}, where $(\xi,\varphi,L)$
satisfies the assumptions (\bf {H.1}), (\bf {H.2}) and (\bf {H.3}). Then, there
exists a positive constant $C_1(\lambda,p,c_0,T)$ such that: $\forall p\in]1,2[$,
\begin{equation}
\begin{array}{ll}\label{estKYZ}
 \E \left[\sup\limits_{t \in [0,T]} |Y_t|^{e^{\lambda t}+1}+ \int_0^T |Z_s|^2 ds + (K_T)^{p}\right] \\ \qquad\quad\leq C_1(\lambda,p,c_0,T)\E \left[1+|\xi|^{e^{\lambda T}+1} +\int_0^T\mid \eta_s \mid^{e^{\lambda s}+1}ds+ \sup\limits_{0\leq t\leq T}\left((L^+_t)^{e^{\lambda t}}\right)^{\frac{p}{p-1}}\right].
 \end{array}
 \end{equation}
\end{proposition}
\noindent
\textbf{Proof.}
We combine (\ref{estY}), (\ref{estZ}) and (\ref{estK}) to get  (\ref{estKYZ}).\qed


\section {Estimate between two solutions}
\subsection{Some useful tools}
We begin with  an estimate for  $\varphi$, whose proof can be found in \cite{BE}.
\begin{lemma}\label{estimatevarphi} If \bf{(H.3)}
 holds  then,
\begin{align*}
&  \E \left[ \dint_0^T   | \varphi(s,Y_s,Z_s)|^{\overline{\alpha}} ds \right] \; \leq \;
C\left(1 + \E\left[\dint_0^T
 {\eta}_s^2 ds\right]+
\E\left[\dint_0^T\vert Z_s\vert^2 ds\right]\right),
\end{align*}
where $\overline{\alpha}=\min(2,
\dfrac{2}{\alpha})$ for $0\leq \alpha<2$, and $C$ is a positive constant which depends on $c_0$ and $T$.
\end{lemma}



\begin{lemma}\label{exist}
{\it There exists a sequence of functions $(\varphi_n)$ such that:
\begin{enumerate}
\item[$(a)$] For each $n$, $\varphi_n$ is bounded and globally Lipschitz
in $(y,z)$ $a.e.$ $t$ and $P$-$a.s.$ $\omega$.
\item[$(b)$] $\displaystyle\sup_{n}\mid \varphi_n(t,\omega, y, z)\mid
\leq |\eta_t|+c_0|z|\sqrt{\vert\ln(|z|)\vert}$, \quad $P$-$a.s.$, $a.e.$ $t\in
[0,T]$.
\item[$(c)$] $\forall N$, $\rho_N (\varphi_n-\varphi)\longrightarrow 0$ as
$n\longrightarrow+\infty$; where $\rho_N(\varphi)=E\left[\int_{0}^{T}\sup\limits_{|y|,|z|\leq N}\mid\varphi(s,y,z)\mid ds\right].$
\end{enumerate}
}
\end{lemma}

\noindent
\textbf{Proof.} Let $\alpha_n: \R^2 \longrightarrow \R_+$ be a sequence of
smooth functions with compact support which approximate the Dirac
measure at 0 and which satisfy $\int \alpha_n (u)du = 1$. Let
$\psi_n$ from $\R^{2}$ to $ \R_+$  be a sequence of smooth functions
such that $0\leq \vert\psi_n\vert \leq 1$, $\psi_n(u)=1$ for $\vert
u\vert \leq n$ and $\psi_n(u)=0$ for $\vert u\vert \geq n+1$.  We
 put, $\varepsilon_{q,n}(t,y,z) = \int
\varphi(t,(y,z)-u)\alpha_q(u)du\psi_n(y,z)$. For $n \in \N^*$, let
$q(n)$ be an integer such that $q(n) \geq n+n^\alpha$. It is not
difficult to see that the sequence $\varphi_n :=
\varepsilon_{q(n),n}$ satisfies all the assertions $(a)$-$(c)$.\qed\\

Using Proposition \ref{estimateK}, Lemma
\ref{estimatevarphi}, Lemma \ref{exist} and standard arguments of
reflected BSDEs, one can prove the following estimates.
\begin{lemma}\label{lem2}
Let $\varphi$, $\xi$ and $L$ be as in Proposition \ref{estimateK}. Let
$(\varphi_n)$ be the sequence of functions associated to $\varphi$
by Lemma \ref{exist}. Denote by $(Y^{n},Z^{n},K^{n})$ the
solution of equation
\begin{equation}\label{imp}
\left\{\begin{array}{l}
Y^{n}_t =\xi + \int_t^T
\varphi_n(s,Y^{n}_s,Z^{n}_s) ds +K^{n}_T-K^{n}_t- \int_t^T Z^{n}_s dB_s,\qquad
\forall t\in[0,T]\\
L^n_t\leq Y^{n}_t,\qquad t\in[0,T]~~~and~~~\int_{0}^{T}\left(Y^{n}_s-L^n_s\right)dK^{n}_s=0,
\end{array}\right.
\end{equation}
where $(L^n_t)_{t\leq T}$ is supposed to be continuous, increasing with respect to $n$ and $\lim\limits_{n\rightarrow+\infty}L^n_t=L_t, \qquad
\forall t\in[0,T].\\$
Then, there exist
constants $\tilde{K}_1$, $\tilde{K}_2$, $\tilde{K}_3$ and $\tilde{K}_4$ such that:
\begin{enumerate}
\item[$a)$] $ \displaystyle\sup_n\E\left[\int_0^T \vert Z_s^{n}\vert^2ds\right]\leq \tilde{K}_1 $.
\item[$b)$] $ \displaystyle \sup_n\E\left[\sup_{0\leq t\leq T}(\mid
Y_t^{n}\mid^{e^{\lambda T}+1})\right] \leq \tilde{K}_2$.
\item[$c)$] $\displaystyle \sup_n\E\left[\dint_0^T \vert \varphi_n(s,Y_{s}^{n},Z_{s}^{n})\vert^{\overline{\alpha}}ds\right] \leq \tilde{K}_{3}$, where $\overline{\alpha}=\min(2,
\dfrac{2}{\alpha})$.
\item[$d)$] $ \displaystyle\sup_n\E \left[(K^{n}_T)^p\right] \leq \tilde{K}_4,\qquad \forall p \in ]1,2[. $
\end{enumerate}
\end{lemma}




\subsection {Estimate between two solutions}
We can now estimate the variation in the solution.


\begin{lemma}\label{lem4}
For every $R\in\N$, $\beta \in
]1,\min\left(3-\frac{2}{\overline\alpha},2\right)[$, $\forall\delta' <
(\beta-1)\min\left(\frac{1}{4M_2^2},\frac{3-\frac{2}{\overline\alpha}-\beta}{2rM_2^2\beta}\right)$\\ and
$\varepsilon>0$, there exists $N_{0}>R$ such that for all $N>N_{0}$ and $T'\leq T$:
\begin{eqnarray}\label{est2s}
&&\limsup_{n,m\rightarrow +\infty} \E \left[\sup_{(T'-\delta')^{+}\leq t\leq T'}\vert
Y_{t}^{n}-Y_{t}^{m}\vert^\beta\right]  + \E\left[
\int_{(T'-\delta')^{+}}^{T'}{\left|
Z_{s}^{n}-Z_{s}^{m}\right|^{2}\over\left(\vert
Y_{s}^{n}-Y_{s}^{m}\vert^{2}+\nu_{R}\right)^{{2-\beta\over
2}}}ds\right]\\ \nonumber
&&\qquad\qquad\qquad\qquad\qquad\qquad\leq \varepsilon +\frac{\ell}{\beta -1}
e^{C_N\delta'} \limsup_{n,m\rightarrow +\infty} \E \left[\vert
Y_{T'}^{n}-Y_{T'}^{m}\vert^\beta\right]
\end{eqnarray}
where $\nu_{R} = \sup \left\{(A_N)^{-1}, N\geq R
\right\}$,
 $C_N = {2M_2^2\beta\over (\beta -1)} \ln A_{N}$, and $\ell$ is a
 universal positive constant.
\end{lemma}

\noindent
\textbf{Proof.} 
Let $0<T'\leq T$. It follows from It\^o's formula that for all
$t\leq T'$,
\begin{align*}
&\left| Y_{t}^{n}-Y_{t}^{m}\right|
^{2}+\int_{t}^{T'}\left| Z_{s}^{n}-Z_{s}^{m}\right|
^{2}ds
\\& = \left| Y_{T'}^{n}-Y_{T'}^{m}\right|
^{2}+2\int_{t}^{T'}\big( Y_{s}^{n}-Y_{s}^{m}\big)
\big(\varphi_{n}(s,Y_{s}^{n},Z_{s}^{n})-
\varphi_{m}(s,Y_{s}^{m},Z_{s}^{m})\big)
ds
\\&+2\int_{t}^{T'}\big( Y_{s}^{n}-Y_{s}^{m}\big)\big(dK^{n}_s-dK^{m}_s \big)-2\int_{t}^{T'}\langle Y_{s}^{n}-Y_{s}^{m},
\quad \left(
Z_{s}^{n}-Z_{s}^{m}\right)dB_{s}\rangle.
\end{align*}
\\
For \;$N\in \N^*$ we set, $ \Delta_{t}:=\left|
Y_{t}^{n}-Y_{t}^{m}\right| ^{2}+ (A_N)^{-1}. $
\\
Let \;$C>0$ and $1<\beta < \min
\{(3-\frac{2}{\overline\alpha}),2\}$. It\^o's formula shows that,
\begin{align*}
& e^{Ct}\Delta_{t}^{\beta\over 2}
+C\int_{t}^{T'}e^{Cs}\Delta_{s}^{\frac{\beta}{2}}ds\\& =
e^{CT'}\Delta_{T'}^{\beta\over 2}
+\beta\int_{t}^{T'}e^{Cs}\Delta_{s}^{\frac{\beta}{2}-1} \big(
Y_{s}^{n}-Y_{s}^{m}\big)
\big(\varphi_{n}(s,Y_{s}^{n},Z_{s}^{n})-
\varphi_{m}(s,Y_{s}^{m},Z_{s}^{m})\big)
ds
\\ & -\frac{\beta}{2}\int_{t}^{T'}e^{Cs}\Delta_{s}^{\frac{\beta}{2}-1}\left|
Z_{s}^{n}-Z_{s}^{m}\right|^{2}ds
-\beta\int_{t}^{T'}e^{Cs}\Delta_{s}^{\frac{\beta}{2}-1}\langle
Y_{s}^{n}-Y_{s}^{m}, \quad \left(
Z_{s}^{n}-Z_{s}^{m}\right)dB_{s}\rangle
\\ &+\beta\int_{t}^{T'}e^{Cs}\Delta_{s}^{\frac{\beta}{2}-1} \big(
Y_{s}^{n}-Y_{s}^{m}\big)\big(dK^{n}_s-dK^{m}_s \big)-\beta(\frac{\beta-2}{2})\int_{t}^{T'}e^{Cs}
\Delta_{s}^{\frac{\beta}{2}-2}
\left((Y_{s}^{n}-Y_{s}^{m})(Z_{s}^{n}-
Z_{s}^{m})\right)^{2}ds.
\end{align*}
Put $\Phi(s)=|Y_{s}^{n}| + |Y_{s}^{m}|+
|Z_{s}^{n}| + |Z_{s}^{m}|+v_s$. Then
\begin{align*}
&
e^{Ct}\Delta_{t}^{\beta\over 2}
+C\int_{t}^{T'}e^{Cs}\Delta_{s}^{\frac{\beta}{2}}ds
\\& =
e^{CT'}\Delta_{T'}^{\beta\over 2}
-\beta\int_{t}^{T'}e^{Cs}\Delta_{s}^{\frac{\beta}{2}-1}\langle
Y_{s}^{n}-Y_{s}^{m}, \quad \left(
Z_{s}^{n}-Z_{s}^{m}\right)dB_{s}\rangle
\\ &
-\frac{\beta}{2}\int_{t}^{T'}e^{Cs}\Delta_{s}^{\frac{\beta}{2}-1}\left|
Z_{s}^{n}-Z_{s}^{m}\right|^{2}ds+\beta\int_{t}^{T'}e^{Cs}\Delta_{s}^{\frac{\beta}{2}-1} \big(
Y_{s}^{n}-Y_{s}^{m}\big)\big(dK^{n}_s-dK^{m}_s \big)
\\ &
+\beta\frac{(2-\beta)}{2}\int_{t}^{T'}e^{Cs}
\Delta_{s}^{\frac{\beta}{2}-2}
\left((Y_{s}^{n}-Y_{s}^{m})(Z_{s}^{n}-
Z_{s}^{m})\right)^{2}ds
\\ &
+J_{1}+J_{2}+J_{3}+J_{4},
\end{align*}
where \begin{align*}
 & J_{1}:=\beta\int_{t}^{T'}e^{Cs}\Delta_{s}^{\frac{\beta}{2}-1}
\big( Y_{s}^{n}-Y_{s}^{m}\big)
\big(\varphi_{n}(s,Y_{s}^{n},Z_{s}^{n})-
\varphi_{m}(s,Y_{s}^{m},Z_{s}^{m})\big)
\1_{\{\Phi(s)>N\}} ds.
 \\ & J_{2}:=\beta\int_{t}^{T'}e^{Cs}\Delta_{s}^{\frac{\beta}{2}-1}
\big( Y_{s}^{n}-Y_{s}^{m}\big)
\big(\varphi_{n}(s,Y_{s}^{n},Z_{s}^{n})-
\varphi(s,Y_{s}^{n},Z_{s}^{n})\big)
\1_{\{\Phi(s)\leq N\}} ds.
\\ & J_{3}:=\beta\int_{t}^{T'}e^{Cs}\Delta_{s}^{\frac{\beta}{2}-1}
\big( Y_{s}^{n}-Y_{s}^{m}\big)
\big(\varphi(s,Y_{s}^{n},Z_{s}^{n})-
\varphi(s,Y_{s}^{m},Z_{s}^{m})\big)
\1_{\{\Phi(s)\leq N\}}ds.
\\ & J_{4}:=\beta\int_{t}^{T'}e^{Cs}\Delta_{s}^{\frac{\beta}{2}-1}
\big( Y_{s}^{n}-Y_{s}^{m}\big)
\big(\varphi(s,Y_{s}^{m},Z_{s}^{m})-
\varphi_{m}(s,Y_{s}^{m},Z_{s}^{m})\big)
\1_{\{\Phi(s)\leq N\}} ds.
\end{align*}
Now we will estimate $J_{1}$, $J_{2}$, $J_{3}$ and $J_{4}$. \\

Let $\kappa
=3-\frac{2}{\overline\alpha}-\beta$. Since $\frac
{(\beta-1)}{2}+\frac{\kappa}{2}+\frac{1}{\overline\alpha}=1$, we
use H\" older's inequality to obtain
\begin{align*}
J_{1} &\leq \beta e^{CT'} \dfrac{1}{N^\kappa}\int_{t}^{T'}
\Delta_{s}^{\frac{\beta-1}{2}}{\Phi^\kappa(s)}
|\varphi_{n}(s,Y_{s}^{n},Z_{s}^{n})-\varphi_{m}(s,Y_{s}^{m},Z_{s}^{m})|ds
\\ &\leq
\beta e^{CT'} \dfrac{1}{N^\kappa} \left[\int_{t}^{T'} \Delta_{s}
ds \right]^{\frac{\beta-1}{2}} \left[\int_{t}^{T'}{\Phi(s)}^2 ds
\right]^{\frac{\kappa}{2}}
\\ &
\times\left[\int_{t}^{T'}|\varphi_{n}(s,Y_{s}^{n},Z_{s}^{n})-
\varphi_{m}(s,Y_{s}^{m},Z_{s}^{m})|^{\overline\alpha}ds
\right]^{\frac{1}{\overline\alpha}}.
\end{align*}
Since $|Y_{s}^{n}-Y_{s}^{m}|\leq
\Delta_s^\frac{1}{2}$, it is easy to see that
\begin{align*}
J_{2}+ J_{4} & \leq 2\beta e^{CT'} [ 2N^2+\nu_1]^{\frac{\beta-1}{2}}
\bigg[\int_{t}^{T'}\sup_{|y|,|z|\leq
N}|\varphi_{n}(s,y,z)-\varphi(s,y,z)|ds
\\  &
+\int_{t}^{T'}\sup_{|y|,|z|\leq
N}|\varphi_{m}(s,y,z)-\varphi(s,y,z)|ds\bigg].
\end{align*}
Using assumption \bf{(H.4)}, we get
\begin{align*}
J_{3} &\leq \beta M_2
\int_{t}^{T'}e^{Cs}\Delta_{s}^{\frac{\beta}{2}-1}
 \bigg[|Y_{s}^{n}-Y_{s}^{m}|^{2}\ln
A_{N}
\\ &
+ \frac{\ln A_N}{A_N} +
|Y_{s}^{n}-Y_{s}^{m}||Z_{s}^{n}-Z_{s}^{m}|\sqrt{\ln
A_{N}} \bigg]\1_{\{\Phi(s)\leq N\}}ds
\\ & \leq
\beta M_2\int_{t}^{T'}e^{Cs}\Delta_{s}^{\frac{\beta}{2}-1}
 \bigg[\Delta_{s}\ln A_{N}+|Y_{s}^{n}-Y_{s}^{m}||Z_{s}^{n}-Z_{s}^{m}|\sqrt{\ln
A_{N}}  \bigg]\1_{\{\Phi(s)\leq N\}}ds.
\end{align*}
Next let us deal with
 $$ \beta\int_{t}^{T'}e^{Cs}\Delta_{s}^{\frac{\beta}{2}-1} \big(
Y_{s}^{n}-Y_{s}^{m}\big)\big(dK^{n}_s-dK^{m}_s \big).$$ Actually, since $dK^n_s=\ind_{\{Y^n_s=L^n_s\}}dK^n_s$ and $dK^m_s=\ind_{\{Y^m_s=L^m_s\}}dK^m_s$ we have
\begin{eqnarray*}
&&\beta\int_{t}^{T'}e^{Cs}\Delta_{s}^{\frac{\beta}{2}-1}\big(
Y_{s}^{n}-Y_{s}^{m}\big)\big(dK^{n}_s-dK^{m}_s \big)\\
&&\qquad\qquad=\beta\int_{t}^{T'}e^{Cs}\left(\left|
Y_{s}^{n}-Y_{s}^{m}\right| ^{2}+ (A_N)^{-1}\right)^{\frac{\beta}{2}-1} \big(
Y_{s}^{n}-Y_{s}^{m}\big)\ind_{\{Y_{s}^{n}=L_{s}^{n}\}}dK^{n}_s\\ \nonumber
&&\qquad\qquad-\beta\int_{t}^{T'}e^{Cs}\left(\left|
Y_{s}^{n}-Y_{s}^{m}\right| ^{2}+ (A_N)^{-1}\right)^{\frac{\beta}{2}-1} \big(
Y_{s}^{n}-Y_{s}^{m}\big)\ind_{\{Y_{s}^{m}=L_{s}^{m}\}}dK^{m}_s.
\end{eqnarray*}
Then, it follows that 
\begin{eqnarray*}
&&\beta\int_{t}^{T'}e^{Cs}\Delta_{s}^{\frac{\beta}{2}-1}\big(
Y_{s}^{n}-Y_{s}^{m}\big)\big(dK^{n}_s-dK^{m}_s \big)\\
&&\qquad\qquad=\beta\int_{t}^{T'}e^{Cs}\left(\left|
L_{s}^{n}-Y_{s}^{m}\right| ^{2}+ (A_N)^{-1}\right)^{\frac{\beta}{2}-1} \big(
L_{s}^{n}-Y_{s}^{m}\big)\ind_{\{L_{s}^{n}-Y_{s}^{m}\neq 0\}}dK^{n}_s\\ \nonumber
&&\qquad\qquad-\beta\int_{t}^{T'}e^{Cs}\left(\left|
Y_{s}^{n}-L_{s}^{m}\right| ^{2}+ (A_N)^{-1}\right)^{\frac{\beta}{2}-1} \big(
Y_{s}^{n}-L_{s}^{m}\big)\ind_{\{Y_{s}^{n}-L_{s}^{m}\neq 0\}}dK^{m}_s.
\end{eqnarray*}
Now, let $F$ be the function $(s,x,y)\mapsto\beta e^{Cs}(|x-y|^2+(A_N)^{-1})^{\frac{\beta}{2}-1} \big(
x-y\big)\ind_{\{x-y\neq 0\}}$.
Therefore, 
\begin{eqnarray*}
&& \beta\int_{t}^{T'}e^{Cs}\Delta_{s}^{\frac{\beta}{2}-1}\big(
Y_{s}^{n}-Y_{s}^{m}\big)\big(dK^{n}_s-dK^{m}_s \big)=\int_{t}^{T'} F(s,L^n_s,Y^m_s)dK^n_s- \int_{t}^{T'} F(s,Y^n_s,L^m_s)dK^m_s.
\end{eqnarray*} 
A simple calculation shows that for every $(x,y)\in \R $ the functions $x\in\R\mapsto F(s,x,y)$ and $y\in\R\mapsto F(s,x,y)$ are respectively non-decreasing and non-increasing. Since $L^m_s\leq Y^m_s$ and $L^n_s\leq Y^n_s$ for every $s\in[0,T]$, we obtain
\begin{eqnarray*}
&&\beta\int_{t}^{T'}e^{Cs}\Delta_{s}^{\frac{\beta}{2}-1} \big(
Y_{s}^{n}-Y_{s}^{m}\big)\big(dK^{n}_s-dK^{m}_s \big)\leq\int_{t}^{T'} F(s,L^n_s,L^m_s)dK^n_s- \int_{t}^{T'} F(s,L^n_s,L^m_s)dK^m_s.
\end{eqnarray*}
It follows that
\begin{multline*}
\beta\int_{t}^{T'}e^{Cs}\Delta_{s}^{\frac{\beta}{2}-1} \big(
Y_{s}^{n}-Y_{s}^{m}\big)\big(dK^{n}_s-dK^{m}_s \big)\\ \leq \beta\int_{t}^{T'}e^{Cs}\left(\left|
L_{s}^{n}-L_{s}^{m}\right| ^{2}+ (A_N)^{-1}\right)^{\frac{\beta}{2}-1} \big(
L_{s}^{n}-L_{s}^{m}\big)\big(dK^{n}_s-dK^{m}_s \big).
\end{multline*}
Put $\Delta L^{n,m}_s=\left(\left|L_{s}^{n}-L_{s}^{m}\right| ^{2}+ (A_N)^{-1}\right)^{\frac{\beta}{2}-1} \big(
L_{s}^{n}-L_{s}^{m}\big),$ then from H\" older's inequality and from Lemma \ref{lem2} there exist $C_1$ and $\tilde{K}_4$ such that: $ \forall p \in ]1,2[$
\begin{eqnarray*}
&&\E\left[\beta\int_{t}^{T'}e^{Cs}\Delta_{s}^{\frac{\beta}{2}-1} \big(
Y_{s}^{n}-Y_{s}^{m}\big)\big(dK^{n}_s-dK^{m}_s \big)\right]\leq \E\left[\beta\int_{t}^{T'}e^{Cs}\Delta L^{n,m}_s\big(dK^{n}_s-dK^{m}_s \big)\right]\\
&&\qquad\qquad\qquad\qquad\qquad\qquad\qquad\qquad\qquad\qquad\leq C_1e^{CT'} \tilde{K}_4^{\frac{1}{p}}\left(\E\left[\sup_{0\leq t\leq T}{|\Delta L^{n,m}_t|^{\frac{p}{p-1}}} \right] \right)^{\frac{p-1}{p}}.
\end{eqnarray*}
We choose $C=C_N=\dfrac{2M_2^2\beta}{\beta -1} \ln A_{N}$, and we use Lemma 4.6 in \cite{BE} and we get:
\begin{eqnarray*}
&& e^{C_Nt}\Delta_{t}^{\beta\over 2} + \dfrac{\beta(\beta-1)}{4}
\int_{t}^{T'}e^{C_Ns}\Delta_{s}^{\frac{\beta}{2}-1}\left|
Z_{s}^{n}-Z_{s}^{m}\right|^{2}ds\\
&&\qquad\qquad\qquad\leq e^{C_NT'}\Delta_{T'}^{\beta\over 2} -\beta\int_{t}^{T'}
e^{C_Ns}\Delta_{s}^{\frac{\beta}{2}-1}\langle
Y_{s}^{n}-Y_{s}^{m}, \quad \left(
Z_{s}^{n}-Z_{s}^{m}\right)dB_{s}\rangle\\
&&\qquad\qquad\qquad+\beta e^{C_NT'} \dfrac{1}{N^\kappa}
\left[\int_{t}^{T'} \Delta_{s} ds \right]^{\frac{\beta-1}{2}}
\times\left[\int_{t}^{T'}{\Phi(s)}^2 ds \right]^{\frac{\kappa}{2}}\\
&&\qquad\qquad\qquad \times\left[\int_{t}^{T'}|\varphi_{n}(s,Y_{s}^{n},Z_{s}^{n})-\varphi_{m}(s,Y_{s}^{m},Z_{s}^{m})|^{\overline\alpha}
ds \right]^{\frac{1}{\overline\alpha}}\\
&&\qquad\qquad\qquad +\beta e^{C_NT'} [2N^2+\nu_1]^{\frac{\beta-1}{2}}
\bigg[\int_{t}^{T'}\sup_{|y|,|z|\leq
N}|\varphi_{n}(s,y,z)-\varphi(s,y,z)|ds\\
&&\qquad\qquad\qquad +\int_{t}^{T'}\sup_{|y|,|z|\leq
N}|\varphi_{m}(s,y,z)-\varphi(s,y,z)|ds\bigg]\\
&&\qquad\qquad\qquad+C_1e^{C_NT'} \tilde{K}_4^{\frac{1}{p}}\left(\E\left[\sup_{0\leq t \leq T}{|\Delta L^{n,m}_t|^{\frac{p}{p-1}}} \right] \right)^{\frac{p-1}{p}}.
\end{eqnarray*}
Burkholder's inequality and H\"older's inequality (since $\frac
{(\beta-1)}{2}+\frac{\kappa}{2}+\frac{1}{\overline{\alpha}}=1$)
allow us to show that there exists a universal constant $\ell>0$
such that $\forall \delta'>0$,

\begin{eqnarray*}
&& \E \left[\sup_{(T'-\delta')^{+}\leq t \leq T'}\left[
e^{C_Nt}\Delta_{t}^{\beta\over 2}\right]\right] + \E\left[
\int_{(T'-\delta')^{+}}^{T'}e^{C_Ns}\Delta_{s}^{\frac{\beta}{2}-1}\left|
Z_{s}^{n}-Z_{s}^{m}\right|^{2}ds\right]\\
&&\qquad\qquad\qquad\leq \frac{\ell}{\beta -1} e^{C_NT'}\bigg\{ \E \left[\Delta_{T'}^{\beta\over 2}\right]
+\dfrac{\beta}{N^\kappa} \left(\E \left[\int_{0}^{T} \Delta_{s} ds
\right]\right)^{\frac{\beta-1}{2}} \left(\left[\E \int_{0}^{T}{\Phi(s)^2} ds\right]
\right)^{\frac{\kappa}{2}}
\\
&&\qquad\qquad\qquad \times\left(\E\left[
\int_{0}^{T}|\varphi_{n}(s,Y_{s}^{n},Z_{s}^{n})-\varphi_{m}(s,Y_{s}^{m},Z_{s}^{m})|^{\overline\alpha}
ds \right]\right)^{\frac{1}{\overline\alpha}}
\\
&&\qquad\qquad\qquad +\beta [2N^2+\nu_1]^{\frac{\beta-1}{2}}
\E \bigg[\int_{0}^{T}\sup_{|y|,|z|\leq
N}|\varphi_{n}(s,y,z)-\varphi(s,y,z)|ds
\\
&&\qquad\qquad\qquad+\int_{0}^{T}\sup_{|y|,|z|\leq
N}|\varphi_{m}(s,y,z)-\varphi(s,y,z)|ds\bigg] \bigg\}
\\
&&\qquad\qquad\qquad+C_1e^{C_NT'} \tilde{K}_4^{\frac{1}{p}}\left(\E\left[\sup_{0\leq t\leq T}{|\Delta L^{n,m}_t|^{\frac{p}{p-1}}} \right] \right)^{\frac{p-1}{p}}.
\end{eqnarray*}
We use Lemma \ref{estimatevarphi}, Lemma \ref{exist} and Lemma \ref{lem2} to obtain, $\forall N>R$,

\begin{align*}
& \E \left[\sup_{(T'-\delta')^{+}\leq t \leq T'}\vert
Y_{t}^{n}-Y_{t}^{m}\vert^\beta\right] + \E\left[
\int_{(T'-\delta')^{+}}^{T'}\dfrac{\left|
Z_{s}^{n}-Z_{s}^{m}\right|^{2}}{\left(\vert
Y_{s}^{n}-Y_{s}^{m}\vert^{2}+ \nu_{R}
\right)^{\frac{2-\beta}{2}}}ds\right]
\\&\leq \frac{\ell}{\beta -1} e^{C_N\delta'} \bigg\{(A_N)^{-\beta\over 2}
+\beta \dfrac{4
\tilde{K}_{3}^{\frac{1}{\overline\alpha}}}{N^{\kappa}}\left(4T\tilde{K}_2
+T\nu_R\right)^{\frac{\beta-1}{2}}\left(8T\tilde{K}_2
+8\tilde{K}_1\right)^{\frac{\kappa}{2}}
\\ & + \E\left[ \vert
Y_{T'}^{n}-Y_{T'}^{m}\vert^\beta \right]+\beta [2N^2+\nu_1
]^{\frac{\beta-1}{2}} \big[\rho_{N}(\varphi_n -
\varphi)+\rho_{N}(\varphi_m - \varphi)\big]\bigg\}\\
& + C_1e^{C_N\delta'} \tilde{K}_4^{\frac{1}{p}}\left(\E\left[\sup_{0\leq t\leq T}{|\Delta L^{n,m}_t|^{\frac{p}{p-1}}} \right] \right)^{\frac{p-1}{p}}.
\\ & \leq
\frac{\ell}{\beta -1} e^{C_N\delta'}\E\left[ \vert
Y_{T'}^{n}-Y_{T'}^{m}\vert^\beta\right]+ \frac{\ell}{\beta
-1}
\dfrac{A_N^{\frac{2M_2^2\delta'\beta}{\beta-1}}}{(A_N)^{\frac{\beta}{2}}}
\\
&+ \frac{4\ell}{\beta -1}\beta
\tilde{K}_{3}^{\frac{1}{\overline\alpha}}\left(4T\tilde{K}_2
+T\nu_R\right)^{\frac{\beta-1}{2}}\left(8T\tilde{K}_2
+8\tilde{K}_1\right)^{\frac{\kappa}{2}}\dfrac{A_N^{\frac{2M_2^2\delta'\beta}{\beta-1}}}
{(A_N)^{\frac{\kappa}{r}}}
\\ & +\frac{\ell}{\beta -1} e^{C_N\delta'}\beta [2N^2+\nu_1
]^{\frac{\beta-1}{2}} \big[\rho_{N}(\varphi_n -
\varphi)+\rho_{N}(\varphi_m - \varphi)\big]
\\&+ C_1e^{C_N\delta'} \tilde{K}_4^{\frac{1}{p}}\left(\E\left[\sup_{0\leq t\leq T}{|\Delta L^{n,m}_t|^{\frac{p}{p-1}}} \right] \right)^{\frac{p-1}{p}}.
\end{align*}
Hence for $ \delta' < (\beta
-1)\min\left(\frac{1}{4M_2^2},\frac{\kappa}{2rM_2^2\beta}\right)$
we derive
\begin{equation*}
\lim_{N\rightarrow+\infty}\dfrac{A_N^{\frac{2M_2^2\delta'\beta}{\beta-1}}}{(A_N)^{\frac{\beta}{2}}}= 0\qquad\qquad \text{and} \qquad\qquad \lim_{N\rightarrow+\infty}\displaystyle\dfrac{A_N^{\frac{2M_2^2\delta'\beta}{\beta-1}}}
{(A_N)^{\frac{\kappa}{r}}}
\displaystyle=0.
\end{equation*}
Also we have $(L^n_t)_{t\leq T}$ is a continuous and increasing sequence, moreover, $\lim\limits_{n\rightarrow+\infty}L_t^n=L_t$. So from Dini's theorem the convergence of $L^n$ is uniform. Next by the dominated convergence, we conclude that
\begin{equation*}
\lim_{(n,m)\rightarrow+\infty}\left(\E\left[\sup_{0\leq s\leq T}{|\Delta L^{n,m}_s|^{\frac{p}{p-1}}} \right] \right)^{\frac{p-1}{p}}=0.
\end{equation*}
Passing to the limits first on $(n,m)$ and next on $N$, and using assertion $(c)$ of Lemma \ref{exist} to obtain the desired result.\qed\\

Now we introduce the comparison theorem.
\begin{theorem}\label{comp}
Let $(\xi,f,L)$ and $(\xi',f',L')$ be two sets of data that satisfies all the assumptions;
(\textbf{H.1}), (\textbf{H.2}), (\textbf{H.3}), and (\textbf{H.4}). And suppose in addition the followings:
\begin{itemize}
\item[(i)] $\xi\leq \xi'$ $P$-a.s.
\item[(ii)] $f(t,y,z)\leq f'(t,y,z)$
$dP\times dt$ a.e.,
$\forall (t,y,z)\in[0,T]\times \R\times \R^d. $
\item[(iii)] $L_t\leq L'_t;\qquad \forall t\in[0,T]$ $P$-a.s.
\end{itemize}
Let $(Y,Z,K)$  be the  solution of the reflected BSDE with data $(\xi,f,L)$ and $(Y',Z',K')$ the solution of the reflected BSDE with data $(\xi',f',L')$. Then,
$$Y_t\leq Y_t',~~~~~~~0\leq t\leq T~~P\text{-a.s.}$$
\end{theorem}

\noindent\textbf{Proof.}
The arguments of this proof are standard. We defer the proof in Appendix.\qed

\begin{remark}\label{rem}
If $L=-\infty$, then $dK=0$ and the comparison theorem holds also in the standard case.
\end{remark}
\section{Existence and uniqueness}
The main result of this section is the following theorem.
\begin{theorem}\label{unique}
  Assume that
\bf{(H.1)}, \bf{(H.2)}, \bf{(H.3)} and \bf{(H.4)} are satisfied. Then, the equation
(\ref{zlogz1}) has a unique solution.
\end{theorem}
\bop \bf {of Theorem \ref{unique}.}
We divide the proof into two steps.\\

\noindent
\textbf{Step 1.} \bf{Existence.}\\

Taking successively $T'=T$,
$T'=(T-\delta')^+$, $T'=(T-2\delta')^{+}...$ in Lemma \ref{lem4}. 
Then we obtain, for every $\beta\in ]1,~
\min\left(3-\dfrac{2}{\overline\alpha}, 2\right)[$
\begin{align}\label{beta1}
\lim_{n,m\rightarrow +\infty}\left( \E \left[\sup_{0\leq t \leq T}\vert
Y_{t}^{n}-Y_{t}^{m}\vert^\beta\right]  + \E\left[
\int_{0}^{T}\dfrac{\left|
Z_{s}^{n}-Z_{s}^{m}\right|^{2}}{\left(\vert
Y_{s}^{n}-Y_{s}^{m}\vert^{2}+\nu_{R}\right)^{\frac{2-\beta}{2}}}ds\right]\right)=
0.
\end{align}
Since $\beta >1$, Lemma \ref{lem2} allows us to show that
\begin{equation}\label{beta}
\lim_{n\rightarrow +\infty} \E\left[\sup_{0\leq t\leq T}\vert
Y_{t}^{n} - Y_{t}\vert^{\beta}\right] = 0.
\end{equation}
Next let us prove that
\begin{equation}\label{exZp}
\lim\limits_{n\rightarrow+\infty}\E\left[\int_{0}^{T}\left| Z_{s}^{n}-Z_{s}^{m}\right|
^{2}ds\right]=0.
\end{equation}
It follows from It\^{o}'s formula that:
\begin{align}\label{Exz}
&\left| Y_{0}^{n}-Y_{0}^{m}\right|
^{2}+\int_{0}^{T}\left| Z_{s}^{n}-Z_{s}^{m}\right|
^{2}ds
\\\nonumber&\qquad = 2\int_{0}^{T}\big( Y_{s}^{n}-Y_{s}^{m}\big)
\big(\varphi_{n}(s,Y_{s}^{n},Z_{s}^{n})-
\varphi_{m}(s,Y_{s}^{m},Z_{s}^{m})\big)
ds
\\\nonumber&\qquad+2\int_{0}^{T}\big( Y_{s}^{n}-Y_{s}^{m}\big)\big(dK^{n}_s-dK^{m}_s \big)-2\int_{0}^{T}\left(Y_{s}^{n}-Y_{s}^{m}\right) \left(
Z_{s}^{n}-Z_{s}^{m}\right)dB_{s}.
\end{align}
First we argue that the third term of the right side in \eqref{Exz} is a martingale. We can deduce from Burkholder-Davis-Gundy's inequality and Lemma \ref{lem2} that there exists a constant $c>0$ such that:
\begin{eqnarray}\label{mar}
&& \E\left[\sup_{0\leq t\leq T}\left|\int_{0}^{t}\left(Y_{s}^{n}-Y_{s}^{m}\right)\left(Z_{s}^{n}-Z_{s}^{m}\right)dB_{s}\right|\right]\\ \nonumber
&&\qquad\qquad\qquad\qquad\qquad\leq c\E \left[\sup_{0\leq s\leq T}|Y_{s}^{n}-Y_{s}^{m}|^2\right]+c\E\left[\int_{0}^{T}|Z_{s}^{n}-Z_{s}^{m}|^2ds\right]\\ \nonumber
&& \qquad\qquad\qquad\qquad\qquad<+\infty.
\end{eqnarray}
Now we deal with the term
$\int_{0}^{T}\big( Y_{s}^{n}-Y_{s}^{m}\big)\big(dK^{n}_s-dK^{m}_s)$. Actually, since $dK^n_s=\ind_{\{Y^n_s=L^n_s\}}dK^n_s$ and $dK^m_s=\ind_{\{Y^m_s=L^m_s\}}dK^m_s$ and since $L^m_s\leq Y^m_s$ and $L^n_s\leq Y^n_s$ for every $s\in[0,T]$, we obtain
\begin{eqnarray}\label{K}
&& \int_{0}^{T}\big( Y_{s}^{n}-Y_{s}^{m}\big)\big(dK^{n}_s-dK^{m}_s)\\\nonumber
&& \qquad\qquad\qquad=\int_{0}^{T}\big( Y_{s}^{n}-Y_{s}^{m}\big)dK^{n}_s-\int_{0}^{T}\big( Y_{s}^{n}-Y_{s}^{m}\big)dK^{m}_s\\\nonumber
&& \qquad\qquad\qquad=\int_{0}^{T}\big( L_{s}^{n}-Y_{s}^{m}\big)\ind_{\{L_s^n\neq Y_s^m\}}dK^{n}_s-\int_{0}^{T}\big( Y_{s}^{n}-L_{s}^{m}\big)\ind_{\{L_s^m\neq Y_s^n\}}dK^{m}_s\\\nonumber
&& \qquad\qquad\qquad\leq\int_{0}^{T}\big( L_{s}^{n}-L_{s}^{m}\big)\ind_{\{L_s^n\neq L_s^m\}}dK^{n}_s-\int_{0}^{T}\big( L_{s}^{n}-L_{s}^{m}\big)\ind_{\{L_s^m\neq L_s^n\}}dK^{m}_s\\\nonumber
&&\qquad\qquad\qquad = \int_{0}^{T}\big( L_{s}^{n}-L_{s}^{m}\big)\big(dK^{n}_s-dK^{m}_s).
\end{eqnarray}
Combining \eqref{Exz}, \eqref{mar} and \eqref{K} to obtain that there exists a constant $c_1$ such that:
\begin{eqnarray}\label{z0}
&& \E \left[\int_{0}^{T}\left| Z_{s}^{n}-Z_{s}^{m}\right|
^{2}ds\right]\\ \nonumber
&& \qquad\qquad\qquad\qquad\leq c_1\E\left[\int_{0}^{T}\big| Y_{s}^{n}-Y_{s}^{m}\big|\big|\varphi_{n}(s,Y_{s}^{n},Z_{s}^{n})-
\varphi_{m}(s,Y_{s}^{m},Z_{s}^{m})\big|ds \right]\\ \nonumber
&& \qquad\qquad\qquad\qquad+c_1\E \left[\int_{0}^{T}\big( L_{s}^{n}-L_{s}^{m}\big)\big(dK^{n}_s-dK^{m}_s)\right].
\end{eqnarray}
Next by  H\" older's inequality we have
\begin{eqnarray}
&& \E\left[\int_{0}^{T}\big| Y_{s}^{n}-Y_{s}^{m}\big|\big|\varphi_{n}(s,Y_{s}^{n},Z_{s}^{n})-
\varphi_{m}(s,Y_{s}^{m},Z_{s}^{m})\big|ds \right]\\ \nonumber
&&\quad\quad\qquad\leq\E\left[\left(\int_{0}^{T}|Y_s^n-Y_s^m|^{\frac{\bar{\alpha}}{\bar{\alpha}-1}}ds\right)^{\frac{\bar{\alpha}-1}{\bar{\alpha}}}\left(\int_{0}^{T}\big|\varphi_{n}(s,Y_{s}^{n},Z_{s}^{n})-
\varphi_{m}(s,Y_{s}^{m},Z_{s}^{m})\big|^{\bar{\alpha}}ds\right)^{\frac{1}{\bar{\alpha}}}\right]\\ \nonumber
&&\quad\quad\qquad\leq\E\left[\int_{0}^{T}|Y_s^n-Y_s^m|^{\frac{\bar{\alpha}}{\bar{\alpha}-1}}ds\right]^{\frac{\bar{\alpha}-1}{\bar{\alpha}}}\E\left[\int_{0}^{T}\big|\varphi_{n}(s,Y_{s}^{n},Z_{s}^{n})-
\varphi_{m}(s,Y_{s}^{m},Z_{s}^{m})\big|^{\bar{\alpha}}ds\right]^{\frac{1}{\bar{\alpha}}}\\ \nonumber
&&\quad\quad\qquad\leq \E \left[\int_{0}^{T}|Y_s^n-Y_s^m|^\beta ds\right]^{\frac{1}{\beta}}\E \left[\int_{0}^{T}|Y_s^n-Y_s^m|^{\frac{\beta}{(\bar{\alpha}-1)(\beta-1)}}ds\right]^{\frac{\beta-1}{\beta}}\\ \nonumber
&&\quad\quad\qquad\times\E\left[\int_{0}^{T}\big|\varphi_{n}(s,Y_{s}^{n},Z_{s}^{n})-
\varphi_{m}(s,Y_{s}^{m},Z_{s}^{m})\big|^{\bar{\alpha}}ds\right]^{\frac{1}{\bar{\alpha}}}.
\end{eqnarray}
We plug the last inequality in \eqref{z0}. Then Proposition \ref{estimateK} (for $\lambda$ large enough), Lemma \ref{estimatevarphi}, Lemma \ref{exist} and Lemma \ref{lem2} shows that for $\beta\in]1,\min\left(3-\frac{2}{\overline{\alpha}},2\right)[$ and $\bar{\alpha}$ large enough there exists a positive constant $C'$ such that:
\begin{eqnarray}
&&\E \left[\int_{0}^{T}\left| Z_{s}^{n}-Z_{s}^{m}\right|
^{2}ds\right]\\ \nonumber
&&\qquad\qquad\leq C'\E\left[\sup_{0\leq t\leq T}|Y_t^n-Y_t^m|^{\beta}\right]^{\frac{\beta-1}{\beta}}+2\tilde{K}_4^{\frac{1}{p}}\E\left[\sup_{0\leq t\leq T}|L_t^n-L_t^m|^{\frac{p}{p-1}}\right]^{\frac{p-1}{p}}.
\end{eqnarray}
From (\ref{beta1})  and the fact that $L^n$ converge uniformly to $L$ it follows that:
\begin{equation}\label{exZp1}
\lim\limits_{n\rightarrow+\infty}\E\left[\int_{0}^{T}\left| Z_{s}^{n}-Z_{s}^{m}\right|
^{2}ds\right]=0.
\end{equation}
Finally we use Lemma \ref{lem2} to show that
\begin{equation}\label{existenceofz}
\lim\limits_{n\rightarrow+\infty}\E \left[\int_{0}^{T}\left| Z_{s}^{n}-Z_{s}\right|
^{2}ds\right]=0.
\end{equation}


On the other hand
\begin{align*}
&\E\left[\dint_{0}^{T}\vert \varphi_{n}(s,Y_{s}^{n},
Z_{s}^{n}) -\varphi(s,Y_{s}^{n},
Z_{s}^{n})\vert ds\right]
\\ & \leq \E\left[\dint_{0}^{T}\vert \varphi_{n}(s,Y_{s}^{n}, Z_{s}^{n})-\varphi(s,Y_{s}^{n},
Z_{s}^{n})\vert \1_{\{\vert Y_{s}^{n}\vert +\vert
Z_{s}^{n}\vert\leq N\}} ds\right]
\\ & +\E\left[\dint_{0}^{T}\vert \varphi_{n}(s,Y_{s}^{n}, Z_{s}^{n})-\varphi(s,Y_{s}^{n},
Z_{s}^{n})\vert\dfrac{(\vert Y_{s}^{n}\vert +\vert
Z_{s}^{n}\vert)^{(2-\frac{2}{\overline\alpha})}}{N^{(2-\frac{2}{\overline\alpha})}}
\1_{\{\vert Y_{s}^{n}\vert +\vert Z_{s}^{n}\vert\geq
N\}} ds\right]
\\ & \leq \rho_{N}(\varphi_n - \varphi)
+\dfrac{2\tilde{K}_{3}^{\frac{1}{\overline\alpha}}\left[T\tilde{K}_{2}
+\tilde{K}_1\right]^{1-\frac{1}{\overline\alpha}}}{N^{(2-\frac{2}{\overline\alpha})}}.
\end{align*}
Passing to the limit first on $n$ and next on $N$  we obtain
\begin{equation}\label{Phi1}
\lim_{n\rightarrow+\infty}\E\left[\dint_{0}^{T}\vert \varphi_{n}(s,Y_{s}^{n},
Z_{s}^{n}) -\varphi(s,Y_{s}^{n},
Z_{s}^{n})\vert ds\right] =0.
\end{equation}
Finally, we use Lemma \ref{exist} and Lemma
\ref{lem2} to show that,
\begin{equation}\label{Phi2}
\lim_{n\rightarrow+\infty} \E\left[\dint_{0}^{T}\vert \varphi_{n}(s,Y_{s}^{n},
Z_{s}^{n}) -\varphi(s,Y_{s}, Z_{s})\vert ds\right] =0.
\end{equation}
For the existence of the process $K$, we have
\begin{eqnarray*}
&& K_t^{{n}}-K_t^{{m}}=
\left(Y_0^{{n}}-Y_0^{{m}}\right)-\left(Y_t^{{n}}-Y_t^{{m}}\right)+\int_{0}^{t}(Z_s^{{n}}-Z_s^{{m}})dB_s\\
&&\qquad\qquad\qquad-\int_{0}^{t}\left(\varphi_{n}(s,Y_{s}^{n},
Z_{s}^{n}) -\varphi_{m}(s,Y_{s}^{m},
Z_{s}^{m})\right)ds.
\end{eqnarray*}
Then
\begin{eqnarray*}
&&\E\left[\sup_{0\leq t\leq T}|K_t^{{n}}-K_t^{{m}}|\right]\leq \E\left[|Y_0^{{n}}-Y_0^{{m}}|+\sup_{0\leq t\leq T}|Y_t^{{n}}-Y_t^{{m}}|+\sup_{0\leq t\leq T}\left|\int_{0}^{t}(Z_s^{{n}}-Z_s^{{m}})dB_s\right|\right]\\
&&\qquad\qquad\qquad\qquad\qquad~~~+\E\left[\int_{0}^{T}|\varphi_{n}(s,Y_{s}^{n},
Z_{s}^{n}) -\varphi(s,Y_{s},
Z_{s})|ds\right]\\
&&\qquad\qquad\qquad\qquad\qquad~~~+\E\left[\int_{0}^{T}|\varphi(s,Y_{s},
Z_{s}) -\varphi_m(s,Y_{s}^{m},
Z_{s}^{m})|ds\right].
\end{eqnarray*}
Consequently from Doob's martingale
inequality, (\ref{beta1}), \eqref{exZp1} and (\ref{Phi2}) we have that:
\begin{equation}\label{KKt}
\lim_{n,m\rightarrow +\infty} \E\left[\sup_{0\leq t\leq T}|K_t^{{n}}-K_t^{{m}}|\right]=0.
\end{equation}
Combining (\ref{KKt}) with Lemma \ref{lem2} allows us to show that:
\begin{equation}\label{KKKt}
\lim_{n\rightarrow +\infty} \E\left[\sup_{0\leq t\leq T}|K_t^{{n}}-K_t|\right]=0.
\end{equation}
Additionally, since $\int_{0}^{T}\left(Y_s^{{n}}-L^n_s\right)dK_s^{{n}}\leq 0$ it follows from (\ref{beta}), (\ref{KKKt}) and $\lim\limits_{n\rightarrow +\infty}L_t^n=L_t$ that $\int_{0}^{T}\left(Y_s-L_s\right)dK_s\leq 0$, and from the fact that $L_t\leq Y_t$ we can deduce that $\int_{0}^{T}\left(Y_s-L_s\right)dK_s=0$.
The existence is proved.
\vskip 0.2cm\noindent

\noindent
\textbf{Step 2.} \bf{Uniqueness.}\\

Let $(Y,Z,K)$ and $(Y',Z',K')$ be two solutions  of
equation (\ref{zlogz1}). 
Arguing as in the proof of Lemma \ref{lem4}, one can show that:
\\
for every $R>2$, $\beta \in
]1,\min\left(3-\dfrac{2}{\overline\alpha},2\right)[$, $\delta' <
(\beta
-1)\min\left(\frac{1}{4M_2^2},\frac{3-\frac{2}{\overline\alpha}-\beta}{2rM_2^2\beta}\right)$
and $\varepsilon>0$\\
 there exists $N_{0}>R$ such that  for all $N>N_{0}$, $\forall T'\leq T$
 \begin{multline*}
  \E \left[\sup_{(T'-\delta')^{+}\leq t \leq T'}\vert
 Y_{t}-Y_{t}^{'}\vert^\beta\right]  + \E\left[
 \int_{(T'-\delta')^{+}}^{T'}\dfrac{\left|
 Z_{s}-Z_{s}^{'}\right|^{2}}{\left(\vert
 Y_{s}-Y_{s}^{'}\vert^{2}+\nu_{R}\right)^{\frac{2-\beta}{2}}}ds\right]\\
\leq
 \varepsilon +\frac{\ell}{\beta -1} e^{C_N\delta'}\E \left[\vert
 Y_{T'}-Y_{T'}^{'}\vert^\beta\right].
\end{multline*}
 Taking successively $T'=T$, $T'=(T-\delta')^+$,
 $T'=(T-2\delta')^{+}...$,  we obtain immediately $Y_t = Y_t '$ and
 also $Z_t = Z_t'$. Finally,
the uniqueness of the process $K$ is deduced from the fact that $Y$ and $Z$ are unique.\qed


\section{Application in mixed
stochastic control with finite horizon}
\medskip

In this section, we use the result on finite horizon reflected BSDEs with one barrier to deal
with the mixed stochastic control problem.
In the sequel $\Omega={\cal C}([0,T],\R^d)$ is the space of
continuous
functions from $[0,T]$ to $\R^d$.

Let us consider a mapping $\sigma:\, (t,\omega)\in[0,T]\times \Omega
\mapsto \sigma(t,\omega)\in\R^{d}\bigotimes \R^{d}$ satisfying the
following assumptions,

\medskip\noindent
$\bf {(1.1)}$ \ \  $\sigma$ is {$\mathcal{P}$}-measurable.

\medskip\noindent
$\bf {(1.2)}$ \ \ There exists a constant $C>0$ such that
$|\sigma(t,\omega)-\sigma(t,\omega')|\leq C||\omega-\omega'||_t$ and $|\sigma(t,\omega)|\leq
C(1+||\omega||_t)$, where for any $(\omega,\omega')\in\Omega^2$ and $t\leq T,\,
||\omega||_t=\sup\limits_{s\leq t}|\omega_s|$.

\medskip\noindent
$\bf{(1.3)}$ \ \ For any $(t,\omega)\in [0,T]\times \Omega$, the matrix
$\sigma(t,\omega)$ is invertible and $|\sigma^{-1}(t,\omega)|\leq
C$ for a positive constant $C$.

\medskip\noindent
Let $x_0\in\R^d$ and $x=(x_t)_{t\leq T}$ be the solution of the
following standard functional differential equation:
\begin{equation}\label{EDS}
x_t=x_0+\int_0^t\sigma(s,x)dB_s,\quad t\leq T;
\end{equation}
 the process $(x_t)_{t\leq T}$ exists, since $\sigma$ satisfies
 $\textbf{(1.1)}-\textbf{(1.3)}$ (see, \cite{RY}, pp. 375). Moreover,
 \begin{equation}\label{estim-eds}
 \E[(||x||_T)^n]<+\infty,\quad \forall n\in[1,+\infty[~~(see,\, \cite{KS},\, \mbox{pp}.\,\,
 306).
 \end{equation}


%

Now let $\mathcal{A}$ be a compact metric space and ${\cal U}$ be the space of
${\cal P}$-measurable processes $u:=(u_t)_{t\leq T}$ with value in
$\mathcal{A}$. Let $f:[0,T]\times \Omega \times \mathcal{A}\mapsto \R^d$ be such
that:
\begin{itemize}
\item[$\bf{(1.4)}$] For each $a\in \mathcal{A}$, the function $(t,\omega)\mapsto f(t,\omega,a)$
is predictable.
\item[$\bf{(1.5)}$] For each $(t,\omega)$, the mapping $a\mapsto f(t,\omega,a)$ is
continuous.
\item[$\bf{(1.6)}$] There exists a real constant $\tilde{K}>0$ such that:
\begin{equation}\label{generateur-cond}
|f(t,\omega,a)|\leq \tilde{K}(1+||\omega||_t),\quad \forall 0\leq t \leq T,\,\omega\in
\Omega,\,a\in \mathcal{A}.
\end{equation}
\end{itemize}
Under the previous assumptions, for a given admissible control strategy $u\in {\cal U}$, the exponential process,

$\displaystyle\Lambda^u_t:= \exp\{\int_0^t\sigma^{-1}(s,x)f(s,x,u_s)dB_s-
\frac{1}{2}\int_0^t|\sigma^{-1}(s,x)f(s,x,u_s)|^2ds\}$, \ \  $0\leq t \leq T$,
\par\noindent
is a martingale. Therefore,
 $\E[\Lambda^u_T]=1$ (see \cite{KS}, pp.
191 and 200 ). The Girsanov theorem guarantees
then that the process
\begin{equation}\label{chang-prob}
B^u_t=B_t-\int_0^t\sigma^{-1}(s,x)f(s,x,u_s)ds,\quad 0\leq t \leq T,
\end{equation}
is a Brownian motion with respect to the filtration ${\cal F}_t$,
under the new probability measure
$$P^u(B)=\E[\Lambda^u_T.\ind_{B}],\quad B\in {\cal F}_T,$$
which is equivalent to $P$. It is now clear from equations
(\ref{EDS}) and (\ref{chang-prob}) that
\begin{equation}\label{eds-contr}
x_t=x_0+\int_0^tf(s,x,u_s)ds+\int_0^t\sigma(s,x)dB^u_s,\quad 0\leq t
\leq T.
\end{equation}
\begin{itemize}
\item[$\bf{(1.7)}$] $h:[0,T]\times \Omega \times \mathcal{A}\mapsto\R$ is
measurable and for each $(t,\omega)$ the mapping $a\mapsto h(t,\omega,a)$
is continuous. In addition, there exists a positive constant $\tilde{K}$ such
that:
\begin{equation}\label{hgenerateur}
|h(t,\omega,a)|\leq \tilde{K}(1+||\omega||_t),\quad \forall 0 \leq t \leq T,\,\omega\in
\Omega,\,a\in \mathcal{A}.
\end{equation}
\item[$\bf{(1.8)}$] $g:[0,T]\times \Omega \mapsto\R$ and $g_1: \Omega \mapsto \R$   are two continuous functions for which there exists a positive constant $C$  such that:
\begin{equation} \label{polycond} |g(t,\omega)|+|g_1(\omega)|\leq
C(1+||\omega||_t),\, \, \forall (t,\omega)\in [0,T]\times \Omega .
\end{equation}
\end{itemize}

We define the payoff $$ J(u,\tau)=\E^u\left[\int_{0}^{\tau}h(s,x,u_s)ds+g(\tau,x)\ind_{\{\tau<T\}}+g_1(x_{\tau})\ind_{\{\tau=T\}}\right],$$
and let us set
\begin{equation}\label{Hamiltonien}H(t,x,z,u_t):= z\sigma^{-1}(t,x)f(t,x,u_t)+h(t,x,u_t)\quad \forall(t,x,z,u_t)\in [0,T]\times\Omega\times\R^{d}\times
\mathcal{A}.\end{equation} The function $H$ is called the Hamiltonian
associated with stochastic control such that :
$ \forall z\in \R^{d}$, the process $(H(t,x,z,u_t))_{t\leq
T}$ is ${\cal P}$-measurable.\\

The Hamiltonian function defined in (\ref{Hamiltonien}) attains its
supremum over the set $\mathcal{A}$ at some $u^* \equiv u^*(t,x,z)\in \mathcal{A}$, for
any given $(t,x,z)\in[0,T]\times \Omega \times \R^d$, namely,
\begin{equation}\label{BN}
\sup\limits_{u\in \mathcal{A}}H(t,x,z,u)=H(t,x,z,u^*(t,x,z)).
\end{equation}
(This is the case, for instance, if the set $\mathcal{A}$ is compact and the
mapping $u\mapsto H(t,x,z,u)$ continuous). Then it can be shown
(see Lemma 1 in Benes \cite{B}), that the mapping $u^*:[0, T]\times
\Omega \times \R^d \mapsto \mathcal{A}$ can be selected to be ${\cal P}
\otimes {\cal
B}(\R^d)$-measurable.\\

Now let $H^*(t,x,z)=\sup\limits_{u\in \mathcal{A}}H(t,x,z,u)$ where $x$ is the
solution of (\ref{EDS}).
\begin{theorem}\label{str}
Let $(Y^*,Z^*,K^*)$ be the solution of the finite horizon reflected BSDE associated with $(g_1(x_T),H^*,g(t,x))$ and $ \tau^*=\inf\{t\in[0,T],Y^*_t<g(t,x))\}\wedge T $,
then $ Y^*_0=J(u^*,\tau^*)$ and $(u^*,\tau^*)$ is the optimal strategy for the controller.
\end{theorem}

\noindent
\textbf{Proof.}
We consider the reflected BSDE associated with $(g_1(x_T),H^*,g(t,x)) $
\begin{equation}\label{Ham}
Y^*_t=g_1(x_T)+\int_{t}^{T}H^*(s,x,Z^*_s)ds+K^*_T-K^*_t-\int_{t}^{T}Z^*_sdB_s.
\end{equation}
Now we will show that the Hamiltonian $H^*$ satisfies  \bf{(H.3)} and \bf{(H.4)}. The proof is actually similar to the one in \cite{BE} but for the sake of the reader we give it again.\\

We begin by showing that $H^*$ satisfies \bf{(H.3)}. Actually it is not difficult to see that $\forall (t,x,z)\in[0,T]\times\Omega\times\R^d$ and $|z|$ large enough, there exist two constant $C>0$ and $c_0>0$ such that:

$$|H^*(t,x,z)|\leq C e^{\|x\|_t}+c_0|z|\sqrt{|\ln(|z|)|}.$$
Next we move on to prove that $H^*$ also satisfies \textbf{(H.4)}. Indeed for every $|y|, |y'|, |z|, |z'|\leq N$ we have:
\begin{eqnarray*}
&&\left(y-y'\right)\left(H^*(t,x,z)-H^*(t,x,z')\right)\ind_{\{v_t(\omega)\leq N\}}\\
&&\qquad\qquad\qquad\qquad\qquad\qquad\qquad\leq|y-y'||H^*(t,x,z)-H^*(t,x,z')|\ind_{\{v_t(\omega)\leq N\}}\\
&&\qquad\qquad\qquad\qquad\qquad\qquad\qquad=|\sigma^{-1}(t,x)||y-y'||z-z'||f(t,x,u^*)|\ind_{\{v_t(\omega)\leq N\}}.
\end{eqnarray*}
Now we take $v_t:=e^{|f(t,x,u^*)|^2}$ and since $\sigma^{-1}$ is bounded by a constant $C$, it follows that:
\begin{eqnarray*}
&&\left(y-y'\right)\left(H^*(t,x,z)-H^*(t,x,z')\right)\ind_{\{v_t(\omega)\leq N\}}\\
&&\qquad\qquad\qquad\qquad\qquad\qquad\leq C|y-y'||z-z'||f(t,x,u^*)|\ind_{ \{ e^{|f(t,x,u^*)|^2} \leq N \} }\\
&&\qquad\qquad\qquad\qquad\qquad\qquad\leq C|y-y'||z-z'|\sqrt{\ln(N)}.
\end{eqnarray*}
It remains to show that $e^{|f(t,x,u^*)|^2}$ belongs to $L^q(\Omega\times[0,T];\R_+)$ for some $q$. Actually, there exists a positive constant $\tilde{K}$ such that:
\begin{eqnarray*}
&& \E\left[\int_{0}^{T}e^{q|f(s,x,u^*)|^2}ds\right]\leq \E\left[\int_{0}^{T}e^{2q\tilde{K}^2(1+\|x\|_s^2)}ds\right]\\
&& \qquad\qquad\qquad\qquad\qquad\leq Te^{2q\tilde{K}^2}\E\left[e^{2q\tilde{K}^2\|x\|^2_T}\right].
\end{eqnarray*}
Since $\sigma$ is with linear growth, the result follows for $q$ small enough.\\

Therefore, by Theorem \ref{unique}, \eqref{Ham} has a unique solution $(Y^*,Z^*,K^*)$. Now since $Y^*_0 $ is a deterministic constant we have,
\begin{eqnarray*}
&& Y^*_0=\E^{u^*}[Y^*_0]=\E^{u^*}\left[g_1(x_T)+\int_{0}^{T}H^*(s,x,Z^*_s)ds+K^*_T-\int_{0}^{T}Z^*_sdB_s\right]\\
&&\qquad\qquad\qquad~=\E^{u^*}\left[Y_{\tau^*}^*+\int_{0}^{\tau^*}H^*(s,x,Z^*_s)ds+K^*_{\tau^*}-\int_{0}^{\tau^*}Z^*_sdB_s\right]\\
&&\qquad\qquad\qquad~=\E^{u^*}\left[Y_{\tau^*}^*+\int_{0}^{\tau^*}h(s,x,u^*)ds+K^*_{\tau^*}-\int_{0}^{\tau^*}Z^*_sdB^{u^*}_s\right].
\end{eqnarray*}
From the definition of $\tau^*$ and the properties of reflected BSDEs we know that
the process $K^*_{\tau^*}$ does not increase between $0$ and $\tau^*$, then $K^*_{\tau^*}=0$. Moreover, $\{\int_{0}^{t}Z^*_sdB^{u^*}_s, t\in [0,T]\}$ is  $P^{u^*}$-martingale then
$$ Y^*_0=\E^{u^*}\left[Y_{\tau^*}^*+\int_{0}^{\tau^*}h(s,x,u^*)ds\right].$$  Since $Y_{\tau^*}^*=g_1(x_T)\ind_{\{\tau^*=T\}}+g(\tau^*,x)\ind_{\{\tau^*<T\}}$, then  we get $Y^*_0=J(u^*,\tau^*)$.
Now, let us consider $u$ an element of ${\cal U}$  and $\tau$  a stopping time. Since $P$ and
$P^{u} $ are equivalent  probabilities on $(\Omega,\F)$ we obtain
\begin{eqnarray*}
&& Y^*_0=\E^{u}[Y^*_0]\\
&&~~~~=\E^{u}\left[Y_{\tau}^*+\int_{0}^{\tau}H^*(s,x,Z^*_s)ds+K^*_{\tau}-\int_{0}^{\tau}Z^*_sdB_s\right]\\
&&~~~~=\E^{u}\left[Y_{\tau}^*+\int_{0}^{\tau}h(s,x,u_s)ds+\int_{0}^{\tau}
\left(H^*(s,x,Z^*_s)-H(s,x,Z^*_s,u_s)\right)ds\right.\\
&&~~~~~~~~~~~~~~\left.+K^*_{\tau}-\int_{0}^{\tau}Z^*_sdB^{u}_s\right].
\end{eqnarray*}
But
$K^*_\tau\geq 0$ and from (\ref{BN}) we have,
$H^*(s,x,Z^*_s)-H(s,x,Z^*_s,u_s)\geq 0$ for any $s\in[0,T].$ On the other hand
\begin{eqnarray*}
&&Y_{\tau}^*=g_1(x_T)\ind_{\{\tau=T\}}+Y_{\tau}^*\ind_{\{\tau<T\}}\\
&&~~~~\geq g_1(x_T)\ind_{\{\tau=T\}}+g(\tau,x)\ind_{\{\tau<T\}}
\end{eqnarray*}
and $ \{\int_{0}^{t}Z^*_sdB^{u}_s, t\in [0,T]\}$ is $P^{u}$-martingale. It follows that
$$J(u^*,\tau^*)=Y^*_0\geq \E^{u}\left[\int_{0}^{\tau}h(s,x,u_s)ds+g(\tau,x)\ind_{\{\tau<T\}}+g_1(x_{T})\ind_{\{\tau=T\}}\right]=J(u,\tau).\qed$$\\

\noindent \textbf{Appendix. Proof of Theorem \ref{comp}}\\ First let $(X_t)_{t\leq T}$ be an $\R$-valued continuous semimartingale and $X_t^+=\max\{X_t,0\}$. From Tanaka's formula we have
\begin{equation*}
dX^+_t=\ind_{\{ X_t>0 \}}dX_t+\frac{1}{2}dL_t,~t\leq T
\end{equation*}
where $(L_t)_{t\leq T}$ is an increasing adapted process such that $\int_{0}^{t}X_sdL_s=0,~\forall t\leq T$.
Next by using It\^o's formula and the fact that $(X^+_t)^2=X_t.X^+_t$ we obtain:
\begin{equation}\label{tan}
d(X^+_t)^2=2X^+_tdX_t+\ind_{\{ X_t>0 \}}d<X,X>_t,~t\leq T.
\end{equation}
Now let $(Y_t,Z_t,K_t)_{t\leq T}$ and $(Y'_t,Z'_t,K'_t)_{t\leq T}$ be the solution of the reflected BSDE with lower barrier associated with $(\xi,f,L)$ and $(\xi',f',L')$ respectively. We put $(\tilde{Y},\tilde{Z},\tilde{K})=(Y-Y',Z-Z',K-K').$\\
For $T'\in[0,T]$, it follows from \eqref{tan} that for all
$t\leq T'$,
\begin{align*}
& |\tilde{Y}_t^+|^{2}
+\int_{t}^{T'}\ind_{\{Y_s> Y'_s\}}| \tilde{Z}_s|
^{2}ds
 =|\tilde{Y}_{T'}^+|^{2}
+2\int_{t}^{T'}\tilde{Y}_s^+
\big(f(s,Y_{s},Z_{s})-
f'(s,Y_{s}',Z_{s}')\big)
ds
\\&\qquad\qquad\qquad\qquad\qquad\qquad~~+2\int_{t}^{T'}\tilde{Y}_s^+d\tilde{K}_s-2\int_{t}^{T'}\tilde{Y}_s^+\tilde{Z}_sdB_{s}.
\end{align*}
We set, $ \Delta_{t}:=|\tilde{Y}_t^+|^{2}+ (A_N)^{-1}\ind_{\{Y_t> Y'_t\}}.$ Then for
$C>0$ and $1<\beta < \min
\{(3-\frac{2}{\overline\alpha}),2\}$, It\^o's formula shows that,
\begin{align*}
& e^{Ct}\Delta_{t}^{\beta\over 2}
+C\int_{t}^{T'}e^{Cs}\Delta_{s}^{\frac{\beta}{2}}ds\\
& =e^{CT'}\Delta_{T'}^{\beta\over 2}
+\beta\int_{t}^{T'}e^{Cs}\Delta_{s}^{\frac{\beta}{2}-1} \tilde{Y}_s^+
\big(f(s,Y_{s},Z_{s})-
f'(s,Y_{s}',Z_{s}')\big)
ds
\\ & -\frac{\beta}{2}\int_{t}^{T'}e^{Cs}\Delta_{s}^{\frac{\beta}{2}-1}\ind_{\{Y_s>Y'_s\}}|\tilde{Z}_s|^{2}ds
-\beta\int_{t}^{T'}e^{Cs}\Delta_{s}^{\frac{\beta}{2}-1}\tilde{Y}_s^+\tilde{Z}_sdB_{s}
\\ &+\beta\int_{t}^{T'}e^{Cs}\Delta_{s}^{\frac{\beta}{2}-1}\tilde{Y}_s^+d\tilde{K}_s-\beta(\frac{\beta-2}{2})\int_{t}^{T'}e^{Cs}
\Delta_{s}^{\frac{\beta}{2}-2}
|\tilde{Y}_s^+|^{2}|\tilde{Z}_s|^{2}ds
\end{align*}
Put $\Phi(s)=|Y_{s}| + |Y_{s}'|+
|Z_{s}| + |Z_{s}'|+v_s$. Then
\begin{align*}
&
e^{Ct}\Delta_{t}^{\beta\over 2}
+C\int_{t}^{T'}e^{Cs}\Delta_{s}^{\frac{\beta}{2}}ds
\\& =
e^{CT'}\Delta_{T'}^{\beta\over 2}
-\beta\int_{t}^{T'}e^{Cs}\Delta_{s}^{\frac{\beta}{2}-1}\tilde{Y}_s^+\tilde{Z}_sdB_{s}
\\ &
-\frac{\beta}{2}\int_{t}^{T'}e^{Cs}\Delta_{s}^{\frac{\beta}{2}-1}\ind_{\{Y_s> Y'_s\}}|\tilde{Z}_s|^{2}ds+\beta\int_{t}^{T'}e^{Cs}\Delta_{s}^{\frac{\beta}{2}-1} \tilde{Y}_s^+d\tilde{K}_s
\\ &
+\beta\frac{(2-\beta)}{2}\int_{t}^{T'}e^{Cs}
\Delta_{s}^{\frac{\beta}{2}-2}
|\tilde{Y}_s^+|^{2}|\tilde{Z}_s|^{2}ds+J_{1}+J_{2}+J_{3},
\end{align*}
where \begin{align*}
 & J_{1}:=\beta\int_{t}^{T'}e^{Cs}\Delta_{s}^{\frac{\beta}{2}-1}
\tilde{Y}_s^+
\big(f(s,Y_{s},Z_{s})-
f'(s,Y_{s}',Z_{s}')\big)
\1_{\{\Phi(s)>N\}} ds.
 \\ & J_{2}:=\beta\int_{t}^{T'}e^{Cs}\Delta_{s}^{\frac{\beta}{2}-1}
\tilde{Y}_s^+
\big(f(s,Y_{s},Z_{s})-
f'(s,Y_{s},Z_{s})\big)
\1_{\{\Phi(s)\leq N\}} ds.
\\ & J_{3}:=\beta\int_{t}^{T'}e^{Cs}\Delta_{s}^{\frac{\beta}{2}-1}
\tilde{Y}_s^+
\big(f'(s,Y_{s},Z_{s})-
f'(s,Y_{s}',Z_{s}')\big)
\1_{\{\Phi(s)\leq N\}}ds.
\end{align*}
Now we will estimate $J_{1}$, $J_{2}$ and $J_{3}$. \\
\\
We start with $J_1$. Let $\kappa
=3-\frac{2}{\overline{\alpha}}-\beta$. Since $\frac
{(\beta-1)}{2}+\frac{\kappa}{2}+\frac{1}{\overline{\alpha}}=1$, we
use H\" older's inequality to obtain
\begin{align*}
J_{1} &\leq \beta e^{CT'} \dfrac{1}{N^\kappa}\int_{t}^{T'}
\Delta_{s}^{\frac{\beta-1}{2}}{\Phi^\kappa(s)}
|f(s,Y_{s},Z_{s})-f'(s,Y_{s}',Z_{s}')|ds
\\ &\leq
\beta e^{CT'} \dfrac{1}{N^\kappa} \left[\int_{t}^{T'} \Delta_{s}
ds \right]^{\frac{\beta-1}{2}} \left[\int_{t}^{T'}{\Phi(s)}^2 ds
\right]^{\frac{\kappa}{2}}
\\ &
\times\left[\int_{t}^{T'}|f(s,Y_{s},Z_{s})-
f'(s,Y_{s}',Z_{s}')|^{\overline{\alpha}}ds
\right]^{\frac{1}{\overline{\alpha}}}.
\end{align*}
\\
Now for $J_2$, since $f(t,y,z)\leq f'(t,y,z)$, we have that $J_{2}\leq 0.$\\
\\
Finally for $J_3$, using assumption \bf{(H.4)}, we get
\begin{align*}
J_{3} &\leq \beta M_2
\int_{t}^{T'}e^{Cs}\Delta_{s}^{\frac{\beta}{2}-1}
 \bigg[|\tilde{Y}_s^+|^{2}\ln
A_{N}+ \frac{\ln A_N}{A_N}\ind_{\{Y_s>Y'_s \}} +
\tilde{Y}_s^+|\tilde{Z}_s|\sqrt{\ln
A_{N}} \bigg]\1_{\{\Phi(s)<N\}}ds
\\ & \leq
\beta M_2\int_{t}^{T'}e^{Cs}\Delta_{s}^{\frac{\beta}{2}-1}
 \bigg[\Delta_{s}\ln A_{N}+\tilde{Y}_s^+|\tilde{Z}_s|\sqrt{\ln
A_{N}}  \bigg]\1_{\{\Phi(s)\leq N\}}ds.
\end{align*}
Next let us deal with
$ \beta\int_{t}^{T'}e^{Cs}\Delta_{s}^{\frac{\beta}{2}-1} \tilde{Y}_s^+\big(dK_s-dK'_s \big)$. Since on $\{Y>Y'\}$, $Y> L'\geq L$ then
 \begin{equation}
\int_{t}^{T'}e^{Cs}\Delta_{s}^{\frac{\beta}{2}-1} \tilde{Y}_s^+\big(dK_s-dK'_s \big)=-\int_{t}^{T'}e^{Cs}\Delta_{s}^{\frac{\beta}{2}-1} \tilde{Y}_s^+dK'_s\leq 0.
 \end{equation}
We apply Lemma 4.6 in \cite{BE}  and we choose $C=C_N=\dfrac{2M_2^2\beta}{\beta -1} \ln A_{N}$, then

\begin{eqnarray*}
&& e^{C_Nt}\Delta_{t}^{\beta\over 2} + \dfrac{\beta(\beta-1)}{4}
\int_{t}^{T'}e^{C_Ns}\Delta_{s}^{\frac{\beta}{2}-1}|\tilde{Z}_s|^{2}\ind_{\{ Y_s>Y'_s \}}ds\\
&&\qquad\qquad\qquad\leq e^{C_NT'}\Delta_{T'}^{\beta\over 2} -\beta\int_{t}^{T'}
e^{C_Ns}\Delta_{s}^{\frac{\beta}{2}-1} \tilde{Y}_s^+\tilde{Z}_sdB_{s}\\
&&\qquad\qquad\qquad+\beta e^{C_NT'} \dfrac{1}{N^\kappa}
\left[\int_{t}^{T'} \Delta_{s} ds \right]^{\frac{\beta-1}{2}}
\times\left[\int_{t}^{T'}{\Phi(s)}^2 ds \right]^{\frac{\kappa}{2}}\\
&&\qquad\qquad\qquad \times\left[\int_{t}^{T'}|f(s,Y_{s},Z_{s})-f'(s,Y_{s}',Z_{s}')|^{\overline{\alpha}}
 ds \right]^{\frac{1}{\overline{\alpha}}}.
\end{eqnarray*}
Burkholder's inequality and H\" older's inequality (since $\frac
{(\beta-1)}{2}+\frac{\kappa}{2}+\frac{1}{\overline{\alpha}}=1$)
allow us to show that there exists a universal constant $\ell'>0$ independent from $N$ such that: $\forall \delta'>0$,

\begin{eqnarray*}
&& \E \left[\sup_{(T'-\delta')^{+}\leq t \leq T'}\left[
e^{C_Nt}\Delta_{t}^{\beta\over 2}\right]\right]\\
&&\qquad\qquad\qquad\leq \ell' e^{C_NT'}\bigg\{ \E \left[\Delta_{T'}^{\beta\over 2}\right]
+\dfrac{\beta}{N^\kappa} \E \left[\int_{0}^{T} \Delta_{s} ds
\right]^{\frac{\beta-1}{2}} \E\left[ \int_{0}^{T}{\Phi(s)^2} ds\right]
^{\frac{\kappa}{2}}
\\
&&\qquad\qquad\qquad \times \E\left[
\int_{0}^{T}|f(s,Y_{s},Z_{s})-f'(s,Y_{s}',Z_{s}')|^{\overline{\alpha}}
ds \right]^{\frac{1}{\overline{\alpha}}}\bigg\}.
\end{eqnarray*}
From Lemma \ref{lem2} there exists an other universal constant $\ell''>0$ that changes from line to line such that: $\forall N>R$,
\begin{align*}
& \E \left[\sup_{(T'-\delta')^{+}\leq t \leq T'}
|\tilde{Y}_t^+|^{\beta}\right]
\leq \ell'' e^{C_N\delta'} \bigg\{(A_N)^{-\beta\over 2}
+ \E\left[
|\tilde{Y}_{T'}^+|^{\beta} \right]+\dfrac{1}{N^\kappa}\bigg\}
\\& \qquad\qquad\qquad\qquad\qquad\leq
\ell'' \left\{e^{C_N\delta'}\E\left[
|\tilde{Y}_{T'}^+|^{\beta}\right]+
\dfrac{A_N^{\frac{2M_2^2\delta'\beta}{\beta-1}}}{(A_N)^{\frac{\beta}{2}}}+\dfrac{A_N^{\frac{2M_2^2\delta'\beta}{\beta-1}}}
{(A_N)^{\frac{\kappa}{r}}}\right\}.
\end{align*}
Hence for $ \delta' < (\beta
-1)\min\left(\frac{1}{4M_2^2},\frac{\kappa}{2rM_2^2\beta}\right)$
we derive
\begin{equation*}
\lim_{N\rightarrow+\infty}\dfrac{A_N^{\frac{2M_2^2\delta'\beta}{\beta-1}}}{(A_N)^{\frac{\beta}{2}}}= 0\qquad\qquad and \qquad\qquad \lim_{N\rightarrow+\infty}\displaystyle\dfrac{A_N^{\frac{2M_2^2\delta'\beta}{\beta-1}}}
{(A_N)^{\frac{\kappa}{r}}}
\displaystyle=0.
\end{equation*}
Then it follows that $\forall \varepsilon>0$ there exists $N_0\in\N$ such that for every $N>N_0$
\begin{equation}\label{eps1}
\E \left[\sup_{(T'-\delta')^{+}\leq t \leq T'}
|\tilde{Y}_t^+|^{\beta}\right]
\leq \ell''e^{C_N\delta'}\E\left[
|\tilde{Y}_{T'}^+|^{\beta}\right]+\varepsilon.
\end{equation}
Now taking successively $T'=T$,
$T'=(T-\delta')^+$, $T'=(T-2\delta')^{+}...$ in \eqref{eps1}, we obtain
\begin{equation}
\E \left[\sup_{0\leq t \leq T}
|(Y_{t}-Y_{t}')^+|^{\beta}\right]=0.
\end{equation}
The proof is now complete.\qed


\end{document}